\documentclass[openany]{book}
\usepackage{amsopn,amstext,amsbsy,amsmath,amscd,amsthm,amsfonts}
\usepackage{color}
\usepackage{url}
\usepackage[english]{babel}
\usepackage{enumitem}
\usepackage{booktabs}
\usepackage{hyperref}
\usepackage{multirow}
\usepackage{afterpage}
\usepackage{float}
\usepackage{graphicx}
\usepackage[export]{adjustbox}

\usepackage[ma,mnf,bachelor]{mccover} 
\usepackage{amsfonts, amsmath, amssymb}
\usepackage{geometry} 

\theoremstyle{plain}
\newtheorem{theorem}                 {\bf Theorem}      [chapter]
\newtheorem{proposition}  [theorem]  {\bf Proposition}

\newtheorem{lemma}        [theorem]  {\bf Lemma}

\theoremstyle{definition}

\newtheorem{definition}   [theorem]  {\bf Definition}

\newtheorem{remark}       [theorem]  {\bf Remark}

\allowdisplaybreaks

\def\nab#1#2#3{\nabla^{\hbox{$\scriptstyle{#1}$}}_{\hbox{$\scriptstyle{#2}$}}{\hbox{$#3$}}}

\def\snab#1#2#3{\hbox{$\nabla$\kern-.1em\raise 1.2 ex\hbox{$\scriptstyle{#1}$}\kern-.5em\lower 0.8 ex\hbox{$#2$}\kern-.0em{$#3$}}}

\def\nsnab#1{\hbox{$\nabla$\kern-.1em\raise 1.0 ex\hbox{$\scriptstyle{#1}$}}}

\def\rn{\mathbb R}

\def\nab#1#2{\hbox{$\nabla$\kern -.3em\lower 1.0 ex
		\hbox{$#1$}\kern -.1 em {$#2$}}}

\def\snab#1#2#3{\hbox{$\nabla$\kern-.1em\raise 1.2 ex\hbox{$\scriptstyle{#1}$}\kern-.5em\lower 0.8 ex\hbox{$#2$}\kern-.0em{$#3$}}}

\def\nsnab#1{\hbox{$\nabla$\kern-.1em\raise 1.0 ex\hbox{$\scriptstyle{#1}$}}}

\def \cn{{\mathbb C}}

\def \rn{{\mathbb R}}

\def \F{\mathcal F}

\def \cn{{\mathbb C}}

\def \rn{{\mathbb R}}

\def \smo{C^{\infty}}

\def\nab#1#2{\hbox{$\nabla$\kern -.3em\lower 1.0 ex
		\hbox{$#1$}\kern -.1 em {$#2$}}}

\def \Re{\mathfrak R\mathfrak e}

\def \ip #1#2{\langle #1,#2 \rangle}
\def \lb#1#2{[#1,#2]}

\def \GLC#1{\mathbf{GL}_{#1}(\cn)}

\def \SLR#1{\mathbf{SL}_{#1}(\rn)}
\def \SL2{\widetilde{\text{\bf SL}}_{2}(\rn)}
\def \slr#1{\mathfrak{sl}_{#1}(\rn)}
\def \SLC#1{\mathbf{SL}_{#1}(\cn)}

\def \SU#1{\text{\bf SU}(#1)}
\def \su#1{\mathfrak{su}(#1)}

\DeclareMathOperator{\trace}{trace}

\allowdisplaybreaks

\def\rn{\mathbb R}

\def \GLC#1{\mathbf{GL}_{#1}(\cn)}

\def \SLR#1{\mathbf{SL}_{#1}(\rn)}
\def \SL2{\widetilde{\text{\bf SL}}_{2}(\rn)}
\def \slr#1{\mathfrak{sl}_{#1}(\rn)}
\def \SLC#1{\mathbf{SL}_{#1}(\cn)}

\def \SU#1{\text{\bf SU}(#1)}
\def \su#1{\mathfrak{su}(#1)}

\begin{document}
\newpage
\thispagestyle{empty}
\maintitle{CMC Tori in the Generalised \\Berger Spheres and their Duals}
\authors{Johanna Marie Gegenfurtner}
\issnum{2022}{8}
\sernum{LUNFMA}{4135}{2022}

\frontcover

\large 
\thispagestyle{empty}
\centerline {\bf\Large Abstract}\vskip1cm
The study of minimal surfaces has a long history, due to the important applications. Given a fixed boundary, one wants to minimise the surface area: this can be used, for example, to minimise the area of the roof of a building. Similarly, looking for constant mean curvature (CMC) provides us with many interesting applications in physics – one of the easiest examples are soap bubbles.
In this work however we occupy ourselves with minimal and constant mean curvature surfaces in the three-dimensional sphere $S^3$ and its dual space $\Sigma^3$.
\smallskip

In Chapter 1 we give a brief overview of the tools of Riemannian and Lorentzian geometry that we will use.
We then take a closer look at $S^3,$ computing its Levi-Civita connection and sectional curvatures: in Chapter 2 with respect to the Riemannian metric $g$ and in Chapter 4 with respect to the Lorentzian metric $h$.
Further, we determine some minimal and CMC tori inside $(S^3,g)$ in Chapter 3 and in $(S^3,h)$ in Chapter 5.\smallskip

We then proceed with the dual space  $\Sigma^3$ of $S^3.$ In Chapter 6, we calculate the Levi-Civita connection and sectional curvatures with respect to $g,$
and with respect to $h$ in Chapter 8.
Again we look for minimal and CMC tori of a certain family in $(\Sigma^3,g)$ in Chapter 7 and in $(\Sigma^3,h)$ in Chapter 9.\smallskip

In the appendix, the reader will find a Maple program. It was written to check the computations of the $S^3$ cases, but it can easily be adapted to $\Sigma^3.$
\vskip6cm
{\it Throughout this work it has been my firm intention to give reference to the stated results and credit to the work of others. All theorems, propositions, lemmas and examples left unmarked are either assumed to be well known, or  are the fruits of my own efforts.}

\newpage 

\newpage 
\thispagestyle{empty}
\centerline {\bf\Large Acknowledgments}
\vskip 1cm
I'd like to thank my supervisor Sigmundur Gudmundsson for his help and guidance. I immensely appreciate the patience and confidence he showed me throughout this cooperation. \smallskip

Further, I'd like to express my deep gratitude towards my family, who has supported me during my studies with relentless encouragement. \smallskip\newline
I'd like to thank my friends for all the beautiful distractions during the work on this project and the previous years.
Lastly, I want to thank {\it la mia amica geniale} Lauren Tropeano for her unfailing warmth and wisdom.
\vskip1pc
\hskip9cm Johanna Marie Gegenfurtner
\phantom{m}


\newpage 
\tableofcontents
\thispagestyle{empty}


\setcounter{page}{0}

\chapter{Introduction}

A main result in Lawson's paper \cite{Lawson} from 1970 is that for any $g\in\mathbb{N},$ there exists an embedding of a  compact minimal surface of genus $g$ into $S^3.$ Moreover, the embedding is not unique if $g$ is not prime. 
Lawson further conjectured that for $g=1,$ the only compact minimal surface embedded in $S^3$ is the flat Clifford torus given by $$T=\{(x_1,x_2,x_3,x_4)\in S^3\subset\rn^4 \ | \ x_1^2+x_2^2=x_3^2+x_4^2=\tfrac{1}{2}\}.$$
This was finally proven by Brendle in 2013, see \cite{Brendle}. It is important to speak of an embedding, since in \cite{Hsiang-Lawson} Hsiang and Lawson  constructed an infinite family of minimally immersed tori into $S^3.$
The Lawson conjecture however fails to hold in the Berger spheres. In \cite{TorralboRot}, Torralbo considered the two-parametric Berger sphere $S^3_b (\kappa,\tau)$ with the metric $$g(A,B)=\frac{4}{\kappa}\cdot(\ip{A}{B}+(\frac{4\tau^2}{\kappa}-1)\ip{A}{X}\ip{B}{X}).$$
Here $X_{(z,w)}=(iz,iw)$ belongs to an orthogonal frame for the tangent bundle $TS^3$ of $S^3$, $\ip{}{}$ denotes the usual metric on the sphere and $\kappa,\tau\in\rn$ with $\kappa>0$ and $\tau\neq0.$ Note that $S^3_b(4,1)$ denotes the standard round sphere.
Torralbo gave an example of another compact embedded minimal torus in the Berger sphere $S^3_b (4, 0.4),$ which is not the Clifford torus $T.$\smallskip

It can also be interesting to study when the mean curvature of a manifold is constant. In \cite{deLima}, De Lima et al. considered a family of tori given by an embedding $\Phi_{r,\tau}:S^1\times S^1\mapsto S_b^3$ into the Berger sphere $S_b^3(4,\tau).$ Here the map $\Phi_{r,\tau}$ is defined by  $$(z,w)\mapsto(z\cdot r, \sqrt{1-r^2}\cdot w),$$ where $r\in(0,1).$
As it turns out each of those tori has constant mean curvature $$H_r=\frac{2r^2-1}{2r\sqrt{1-r^2}}.$$\smallskip

In this work we will consider a family of tori parametrised analogously to the above, we however generalise the Berger metric further and obtain a three parameter family of metrics. Then we check under which conditions this matches the result of De Lima et al. in \cite{deLima}. Then we define a Lorentzian Berger metric and investigate the mean curvature. Additionally we will look for constant mean curvature tori of a certain family in $$\Sigma^3=\{(z,w)\ | \ z,w\in\cn^2,\ |z|^2-|w|^2=1\},$$
under both the Riemannian and Lorentzian Berger metrics.\smallskip

As it will be shown, the switch from the Riemannian to Lorentzian metrics does not impact the mean curvature H of the tori in question inside $S^3$ and $\Sigma^3$ respectively.\smallskip

Further we show that under certain conditions, for every real $C\geq0,$ there exist tori in $S^3$ with constant mean curvature $$\|H_g\|\equiv C.$$ For the tori in $\Sigma^3,$ this can only be shown for $C>\tfrac{1}{\mu},$ where $\mu$ is a parameter belonging to the generalised Berger metric. A similar statement, but using a different generalisation of the Berger metric, was given by Torralbo in \cite{TorralboComp}.\smallskip

We also compute the Levi-Civita connection and sectional curvatures in all four cases.

\chapter{Basics on Riemannian/Lorentzian Geometry}
We assume that the reader is familiar with fundamental Riemannian geometry, however we first review some essential facts, that will be needed in our study of minimal and CMC submanifolds. This chapter is largely based on do Carmo's introductory textbook \cite{doCarmo} and Gudmundsson's lecture notes \cite{Gud-Rie}, whereas O'Neill's textbook \cite{O'Neill} was used for the parts on Lorentzian and semi-Riemannian geometry.\smallskip

A {\it differentiable manifold} $(M^m,\mathcal{A})$ is defined to be a topological manifold $M$ of dimension $m$ together with a family of local charts $\mathcal{A}.$ We require $M$ to be locally homeomorphic to $\mathbb R^m,$ i.e. for every $p\in M$ there exists a neighbourhood $U_p$ containing $p$ and a homeomorphism $$x_p: U_p\mapsto \mathbb R^m.$$ Further, the atlas $$\mathcal{A}=\{(U_\alpha,x_\alpha), \alpha\in\mathcal{I}\},$$ needs to cover the whole of $M$ and must be maximal, subject to the condition that the transition maps $$x_\beta \circ x_\alpha^{-1}\vert_{x_\alpha(U_\alpha\cap U_\beta)}:x_\alpha(U_\alpha\cap U_\beta)\subset\rn^m\to\rn^m, \quad \forall \ \alpha,\beta\in\mathcal{I}$$ are differentiable of class $C^\infty.$\smallskip

Let $(M^m,\mathcal{A}_M)$ be a differentiable manifold, where $\mathcal{A}_M$ is an atlas on $M.$ Following Proposition 2.11 in \cite{Gud-Rie}, a subset $N$ of $M$, is called a {\it submanifold} of $M$ if for each $p\in N,$ $\mathcal{A}_M$ contains some $(U_p,x_p),$ such that $p\in U_p$ and $$x_p(U_p\cap N)=x_p(U_p)\cap(\rn^n\times\{0\}).$$ Here $n\in\mathbb{N}$ such that $n\leq m$ denotes the dimension of $N,$ and $m-n$ is called the codimension of $N$ in $M.$ The atlas $\mathcal{A}_M$ induces a structure on N, denoted $\mathcal{A}_N,$ see Proposition 2.11 in \cite{Gud-Rie}. \smallskip

Let $\{e_1,\dots,e_n\}$ be an orthonormal basis of a semi-Euclidean space $V,$ and for a nondegenerate scalar product $\ip{}{},$ set 
$$\epsilon_i=\ip{e_i}{e_i}=\pm 1.$$
The number $v,$ where $0\leq v\leq n,$ of negative signs in the {\it signature} $(\epsilon_1,\dots,\epsilon_n)$ is called the {\it index} of $V.$ Note that the index is independent of the choice of basis, as shown in Lemma 26 in \cite{O'Neill}.\smallskip\newline

Let $M$ be a smooth manifold. A metric tensor $h$ is a symmetric and nondegenerate tensor field of type (0,2) of constant index. That is, to every point $p\in M,$ $h$ associates a scalar product $$h_p: T_p M\times T_p M\mapsto \mathbb R,$$
such that the index of $h_p$ is the same for all $p.$ The pair $(M,h)$ is a 
{\it semi-Riemannian manifold}. A metric on a submanifold $N\subset M$ is obtained by restricting $h$ to $N.$ (Also see Definition 5.5 in \cite{Gud-Rie}.)
\smallskip 

A {\it Riemannian metric} has index 0 and is positive-definite. If $v=1$ and $\dim M\geq2,$ $M$ is a so called {\it Lorentz manifold}. 
Particularly in the context of relativity we might want to describe the {\it causal character} of a tangent vector $W$ of a semi-Riemannian manifold $(M,h):$ 
If 
\begin{eqnarray*}
h(W,W)<0,& &\textrm{$W$ is {\it timelike,}}\\
h(W,W)=0,& \ \textrm{and $W\neq0$,}& \textrm{$W$ is {\it lightlike} }\\
h(W,W)>0,& \ \textrm{or $W=0$},& \textrm{$W$ is  {\it spacelike}}.
\end{eqnarray*}
According to Lemma 25 in \cite{O'Neill} each vector $W\in V$ can uniquely be written as $$W=\sum_{i=1}^{n} \epsilon_i\cdot h(W,e_i)\cdot e_i.$$
To distinguish the Riemannian from the Lorentzian case, we will write $g$ for a Riemannian metric and $h$ for its Lorentzian counterpart.\smallskip

For a differentiable manifold $M,$ and vector fields $X,Y\in\smo(TM),$ and $p\in M,$ the {\it Lie bracket} $$[X,Y]_p:\smo(M)\to\rn$$ is defined by $$X_p(Y(f))-Y_p(X(f)).$$\smallskip

We recall that the {\it Levi-Civita connection} $$\nabla: C^\infty(TM)\times C^\infty(TM)\mapsto C^\infty(TM)$$ on a Riemannian manifold $(M,g)$  is defined to be the unique affine connection that is symmetric and compatible with the Riemannian metric. Given an orthonormal frame $\{E_1,\dots,E_n\}$ of the tangent bundle, the Levi-Civita connection for $A,B\in C^\infty (TM)$ is given by $$\nab AB=\sum_{i=1}^{n} g(\nab AB,E_i)\cdot E_i.$$
The coefficients are given by the {\it Koszul formula}: For $A, B, C\in C^\infty (TM),$
\begin{eqnarray}\label{Koszul}
	g(\nab AB,C)&=&\frac{1}{2}\cdot\{
	A(g(B,C))+B(g(C,A))-C(g(A,B))\\
	&&\quad +g(C,\lb AB)+g(B,\lb CA)-g(A,\lb BC)\}.\nonumber
\end{eqnarray}
We will later use the Levi-Civita connection to compute the {\it sectional curvatures} of $S^3$ and $\Sigma^3.$ Sectional curvatures describe the curvature of Riemannian manifold of dimension larger than 1. The {\it Riemann curvature operator} $R$ is given by \begin{eqnarray}\label{RCO}
R(A,B)C&=&\nab {A}{\nab BC}-\nab {B}{\nab AC}-\nab{\lb AB}{C}.\end{eqnarray}
The sectional curvature at a point $p$ is given by \begin{equation}\label{sectionalcurvature}K_p(A,B)=\frac{g(R(A,B)B,A)}{|A|^2|B|^2-g(A,B)^2}.\end{equation}
\smallskip

We now define two operators on a submanifold $N$ of $M.$ For this, first observe the following:
Let $X$ be a smooth vector field on M, and $\tilde X$ be the restriction of $X$ to the submanifold $N.$ Then for each $p\in N,$ the tangent vector $\tilde X_p\in T_pM$ can uniquely be split into $$\tilde X_p=\tilde X_p^\top  +\tilde X_p^\perp,$$ where $\tilde X_p^\top \in T_pN$ and $\tilde X_p^\perp\in N_pN.$\smallskip

Let now $N$ be a submanifold of $(M,g)$ with induced metric, and let $X,Y$ be local extensions of $\tilde X, \tilde Y\in C^\infty (TN)$ to $\smo(TM).$
The Levi-Civita connection $\tilde \nabla$ on the submanifold $N$ is the part of $\nab XY$ lying in the tangent space of M, i.e. $$\tilde \nabla_{\tilde X}{\tilde Y}=(\nab XY)^\top ,$$ as stated in Definition 6.20 in \cite{Gud-Rie}.

Further, we define the {\it second fundamental form} of $N$ in $M,$ $$B: C^\infty(TN)\otimes C^\infty(TN)\mapsto C^\infty (NN),$$ given by 
$$B(\tilde X, \tilde Y)=(\nab XY)^\perp.$$

For an orthonormal basis $\{E_1,\dots,E_n\}$ of $T_pN,$ the mean curvature vector $H$ is given by $$H=\tfrac{1}{n}\cdot\trace B=\tfrac{1}{n}\sum^n_{i=1} B(E_i,E_i).$$ 
We say that a manifold is minimal if $\trace B=0$ holds everywhere.\smallskip


\chapter{The Generalised Riemannian Berger Spheres $(S^3,g)$}

In this chapter we introduce the 3-dimensional unit sphere $S^3\subset \mathbb R^4\cong\mathbb C^2$ and the generalised Berger metric $g.$ We define a multiplication $\cdot$ on $S^3$ and show that $(S^3,\cdot)$ is isomorphic to the Lie group $(\SU 2,*).$ Using the properties of matrix Lie groups, we compute the Levi-Civita connection as well as the sectional curvatures. \smallskip

Equip the complex two dimensional vector space $\cn^2$ with the standard scalar product on $\rn^4,$ $$\ip{(z_1,w_1)}{(z_2,w_2)} = \Re (z_1\bar z_2 + w_1\bar w_2)$$ and consider the 3-dimensional unit sphere $S^3$ in $\cn^{2}$ given by $$S^3=\{(z,w)\in\cn^2\,|\, |z|^2+|w|^2=1\}.$$
We show that $S^3$ is isomorphic to the Lie group $\SU 2.$
The map $ \Phi:S^3\to \GLC 2$, given by $$
\Phi:(z,w)\mapsto 
\begin{pmatrix}
z & -\bar w\\ 
w & \bar z
\end{pmatrix}.
$$ is an embedding of $S^3$ into the set of invertible $2\times 2$ matrices. We see that the image of $S^3$ under $\Phi$ is the special unitary group $\SU2$ in $\GLC 2$, defined as follows: 
$$\SU2=\{\begin{pmatrix}
z & -\bar w\\ 
w & \bar z
\end{pmatrix} \,|\, |z|^2+|w|^2=1, \ z,w\in\cn^2\}. $$

Let $*$ denote the standard matrix multiplication on $\SU2.$
Then
\begin{eqnarray*}
\Phi(z_1,w_1)*\Phi(z_2,w_2) &=& \begin{pmatrix}z_1 & -\bar w_1\\ w_1 & \bar z_1\end{pmatrix}*\begin{pmatrix}
z_2 & -\bar w_2\\ 
w_2 & \bar z_2
\end{pmatrix} \\\\
&=& \begin{pmatrix}
    z_1z_2-\bar w_1w_2 & -z_1\bar w_2-\bar w_1\bar z_2\\
    \bar z_1w_2+w_1z_2 & \bar z_1\bar z_2-w_1\bar w_2 
\end{pmatrix}
\\ \\
&=&\Phi(z_1z_2-\bar w_1w_2, \bar z_1w_2+w_1z_2).
\end{eqnarray*}
This induces a group structure $\cdot$ on $S^3.$ In particular, $$(z_1,w_1)\cdot (z_2,w_2)=(z_1z_2-\bar w_1w_2, \bar z_1w_2+w_1z_2). $$
Consequently $ \Phi:(S^3,\cdot)\to (\SU 2,*)$ is a group isomorphism. This turns $(S^3, \cdot)$ into a Lie group with neutral element $e=(1,0).$ The multiplicative inverse of an element $(z,w)$ is given by $(\bar z,-w).$\smallskip

The Riemannian metric $\ip{}{}$ on the vector space $\cn^{n\times n}$ of complex $n\times n$ matrices is given by $$\ip XY=\tfrac 12\cdot\mathfrak{Re}(\trace(\bar X^t\cdot Y)),$$
where $X,Y\in\cn^{n\times n}.$ \smallskip

At the neutral element $e=(1,0),$ $T_e\rn^4$ has the orthogonal decomposition 
$$T_e\rn^4=T_eS^3\oplus N_eS^3.$$
First we want to determine a basis for the tangent space of $T_eS^3\cong T_e\SU 2$. By Theorem 3.13 in \cite{Gud-Rie}, $$T_e\SU 2=\{X\in\cn^{2\times2} \,|\, \trace X=0, \ \bar X^t+X=0\}.$$
From this it can be shown that an orthonormal basis for the tangent space $T_e S^3$ of $S^3$ at the unit element $e$ is given by $(i,0),(0,-1),(0,i).$
The normal space $N_e S^3$ is spanned by $N_e=(1,0).$
\begin{definition}\label{metricg}
We now introduce the family of left-invariant Riemannian metrics $$\{g:\smo (TS^3)\otimes 
\smo (TS^3)\mapsto \smo(S^3)\,|\,\lambda \ ,\mu \ ,\nu\in\rn^+\},$$ associating each $p\in S^3$ with a real scalar product
$$g_p:T_p S^3\otimes 
T_p S^3\to\rn$$ 
defined by 
\begin{eqnarray*}
&&g_p(A,B)\\
&=&\lambda^2\cdot\ip{p^{-1} A}{(i,0)}\ip{p^{-1} B}{(i,0)}+\mu^2\cdot\ip{p^{-1}A}{(0,-1)}\ip{p^{-1} B}{(0,-1)}\\
& &\quad +\,\nu^2\cdot\ip{p^{-1} A}{(0,i)}\ip{p^{-1} B}{(0,i)}+\ip{p^{-1} A}{(1,0)}\ip{p^{-1} B}{(1,0)},
\end{eqnarray*}
where $A,B\in T_pS^3.$
\end{definition}
\begin{remark}
The classic Berger metric is obtained by setting $\mu=\nu=1.$ The standard metric is obtained by additionally setting $\lambda=1.$
\end{remark}
Also note that for $A,B\in T_pS^3,$ the last term $\ip{p^{-1} A}{(1,0)}\ip{p^{-1} B}{(1,0)}=0,$ since $N_e$ is normal to the tangent space $T_eS^3$. However, to avoid confusion, we want to define $g$ as a metric on $C^\infty (T\mathbb{R}^4),$ so that we can use it throughout the thesis.

\begin{definition}\label{basis}
With respect to the generalised Berger metric $g$ as defined above, we define the vector fields $X,Y,Z\in\smo(TS^3)$ on $S^3$ by
$$
X_p=\lambda^{-1}\cdot p\cdot(i,0), \ Y_p=\mu^{-1}\cdot p\cdot(0,-1), \ Z_p=\nu^{-1}\cdot p\cdot(0,i),
$$
such that $\{X_p, Y_p, Z_p\}$ forms an orthonormal frame for the tangent bundle $T S^3.$
The normal bundle $N S^3$ is spanned by the vector field $N_p=p\cdot(1,0).$
\end{definition}
Using standard matrix multiplication on $\su 2$ we obtain the Lie bracket relations
$$\lb{X}{Y}=2\,\lambda^{-1}\mu^{-1}\nu\cdot Z,\quad \lb{Z}{X}=2\,\nu^{-1}\lambda^{-1}\mu\cdot Y,\quad \lb{Y}{Z}=2\,\mu^{-1}\nu^{-1}\lambda\cdot X.$$
This follows from basic properties of the Lie bracket. For further clarification, we refer the reader to Chapter 4 of \cite{Gud-Rie}, in particular Propositions 4.33 and 4.36.
The Levi-Civita connection $\nabla$ for $(S^3,g)$ can now be calculated using the Koszul formula \ref{Koszul}


\begin{lemma}\label{LemmaLC}
Equip $S^3$ with the generalised Berger metric $g,$ as given above. Then the Levi-Civita connection $\nabla$ satifies
$$\nab XX=0,\quad \nab XY=\frac{-\lambda^2+\mu^2+\nu^2}{\lambda\mu\nu}\cdot Z,\quad \nab XZ=-\frac{-\lambda^2+\mu^2+\nu^2}{\lambda\mu\nu}\cdot Y,$$
\smallskip$$\nab YX=\frac{-\lambda^2+\mu^2-\nu^2}{\lambda\mu\nu}\cdot Z,\quad \nab YY=0,\quad \nab YZ=-\frac{-\lambda^2+\mu^2-\nu^2}{\lambda\mu\nu}\cdot X,$$
\smallskip$$\nab ZX=\frac{\lambda^2+\mu^2-\nu^2}{\lambda\mu\nu}\cdot Y,\quad \nab ZY=-\frac{\lambda^2+\mu^2-\nu^2}{\lambda\mu\nu}\cdot X,\quad \nab ZZ=0.$$
\end{lemma}

\begin{proof}
Since $X,Y,Z$ is an orthonormal frame for the tangent bundle, we have \begin{equation}\label{nabAB}\nab{A}{B}=g(\nab{A}{B},X)X+g(\nab{A}{B},Y)Y+g(\nab{A}{B},Z)Z.\end{equation}
Note that $g$ is a left-invariant metric and $X,Y,Z$ are left-invariant vector fields, hence by Proposition 6.13 in \cite{Gud-Rie} the entire first row of the Koszul formula vanishes, which yields 	$$2\cdot g(\nab AB,C)=g(C,\lb AB)+g(B,\lb CA)-g(A,\lb BC).$$
This means that $\nab AB$ is entirely determined by the Lie bracket $[,]$ and the Riemannian metric $g$. 
In the following, we explicitly compute $\nab XY,$ the other cases are similar.
\begin{eqnarray*}
2\cdot g(\nab XY,X)&=&g(X,\lb XY)+g(Y, \lb XX)-g(X,\lb YX)\\
&=&g(X,2\,\lambda^{-1}\mu^{-1}\nu\cdot Z)+g(Y,0)-g(X,-2\,\lambda^{-1}\mu^{-1}\nu\cdot Z)\\
&=&0,\\
2\cdot g(\nab XY,Y)&=&g(Y,\lb XY)+g(Y,\lb YX)-g(X,\lb YY)\\&=&g(Y,2\,\lambda^{-1}\mu^{-1}\nu\cdot Z)+g(Y,-2\,\lambda^{-1}\mu^{-1}\nu\cdot Z)-g(X,0)\\
&=&0,\\
2\cdot g(\nab XY,Z)&=&g(Z,\lb XY)+g(Y,\lb ZX)-g(X,\lb YZ)\\
&=&g(Z,2\,\lambda^{-1}\mu^{-1}\nu\cdot Z)+g(Y,2\,\nu^{-1}\lambda^{-1}\mu\cdot Y)-g(X,2\,\mu^{-1}\nu^{-1}\lambda\cdot X)\\
&=&2\,\lambda^{-1}\mu^{-1}\nu+2\,\nu^{-1}\lambda^{-1}\mu-2\,\mu^{-1}\nu^{-1}\lambda\\
&=&\frac{2}{\lambda\mu\nu}\cdot(-\lambda^2+\mu^2+\nu^2).
\end{eqnarray*}
By plugging the above into equation (\ref{nabAB}), we obtain $$\nab XY=g(\nab XY,Z)\cdot Z=\frac{-\lambda^2+\mu^2+\nu^2}{\lambda\mu\nu}\cdot Z.$$
\end{proof}

\begin{proposition}\label{prop:curvatures3}
Let $(S^3,g)$ be a generalised Berger Sphere. Then the sectional curvatures satisfy 
\begin{eqnarray*}
K(X,Y)&=&g(R(X,Y)Y,X)=\frac{(\lambda^2-\mu^2+\nu^2)^2+4\nu^2(\mu^2-\nu^2)}{(\lambda\mu\nu)^2},\\
K(X,Z)&=&g(R(X,Z)Z,X)=\frac{(\lambda^2+\mu^2-\nu^2)^2-4\mu^2(\mu^2-\nu^2)}{(\lambda\mu\nu)^2},\\ 
K(Y,Z)&=&g(R(Y,Z)Z,Y)=\frac{(\lambda^2+\mu^2+\nu^2)^2-4(\lambda^4+\mu^2\nu^2)}{(\lambda\mu\nu)^2}.
\end{eqnarray*}
\end{proposition}
\begin{proof}
We will only prove the first identity, the other ones are computed similarly. For this, we determine the Riemann curvature operator as given in (\ref{RCO}): 
 \begin{eqnarray*}
R(X,Y)Y&=&\nab {X}{\nab YY}-\nab {Y}{\nab XY}-\nab{\lb XY}{Y}\\
&=&-\frac{-\lambda^2+\mu^2+\nu^2}{\lambda\mu\nu}\cdot\nab YZ-\frac{2\nu^2}{\lambda\mu\nu}\cdot\nab ZY\\
&=&\frac{-\lambda^2+\mu^2+\nu^2}{\lambda\mu\nu}\cdot\frac{-\lambda^2+\mu^2-\nu^2}{\lambda\mu\nu}\cdot X+\frac{2\nu^2}{\lambda\mu\nu}\cdot\frac{\lambda^2+\mu^2-\nu^2}{\lambda\mu\nu}\cdot X\\
&=&\frac{\lambda^4-2\lambda^2\mu^2+\mu^4-3\nu^4+2\mu^2\nu^2+2\lambda^2\nu^2}{(\lambda\mu\nu)^2}\cdot X\\
&=&\frac{(\lambda^2-\mu^2+\nu^2)^2+4\nu^2(\mu^2-\nu^2)}{(\lambda\mu\nu)^2}\cdot X.
\end{eqnarray*}
For the values of the Levi-Civita connection we used the results of Lemma \ref{LemmaLC}. Note that the denominator of the sectional curvature as given in (\ref{sectionalcurvature}) simplifies to 1, since we chose an orthonormal basis. Finally, $$K(X,Y)=g(R(X,Y)Y,X)=\frac{(\lambda^2-\mu^2+\nu^2)^2+4\nu^2(\mu^2-\nu^2)}{(\lambda\mu\nu)^2}.$$
\end{proof}
\begin{remark}
 Unsurprisingly, if we set $\lambda=\mu=\nu=1,$ we obtain the round sphere with respect to its standard metric, which has constant sectional curvature 1.\smallskip
 
Generally, 
\begin{equation}\label{eq:KXY31}
    K(X,Y)\leq0
\end{equation} if $\mu<\nu$ and
\begin{equation}\label{eq:KXY32}
    0<\lambda^2\leq 2\nu\sqrt{\nu^2-\mu^2}+\mu^2-\nu^2,
\end{equation}
with equality in (\ref{eq:KXY31}) only if we have equality in  (\ref{eq:KXY32}). 
\begin{equation}\label{eq:KXZ31}
    K(X,Z)\leq0
\end{equation}
if $\mu>\nu$ and
\begin{equation}\label{eq:KXZ32}
    0<\lambda^2\leq 2\mu\sqrt{\mu^2-\nu^2}+\nu^2-\mu^2,
\end{equation}
with equality in (\ref{eq:KXZ31}) only if we have equality in  (\ref{eq:KXZ32}). 
\begin{equation}\label{eq:KYZ31}
    K(Y,Z)\leq0
\end{equation} if \begin{equation}\label{eq:KYZ32}
    \lambda^2\geq\frac{1}{3}\cdot(2\sqrt{\mu^4-\mu^2\nu^2+\nu^4}+\mu^2+\nu^2),
\end{equation}
with equality in (\ref{eq:KYZ31}) only if we have equality in  (\ref{eq:KYZ32}). 
\end{remark}
\chapter{CMC Tori in the Riemannian $(S^3,g)$}
In this chapter we look at a family of tori inside the unit sphere $S^3$ and investigate for which parameters they have a constant mean curvature, and in particular, when they are minimal. We then compare our results to some of the papers mentioned in the introduction. \smallskip\newline

For $\theta,\alpha,\beta\in\rn$, the three dimensional unit sphere $S^3$ can be expressed as $$S^3=(\cos\theta\cdot e^{i\alpha},\sin\theta\cdot e^{i\beta}).$$
If we now fix $0<\theta<\frac{\pi}{2}$, we obtain a two dimensional torus $T^2_\theta$, parametrised by $\mathcal{F}_\theta:\rn^2\to S^3,$ given by $$\F_\theta: (\alpha,\beta)\to(\cos\theta\cdot e^{i\alpha},\sin\theta\cdot e^{i\beta}).$$
It can easily be seen that setting $\theta=0$ or $\theta=\frac{\pi}{2}$ is rather uninteresting, since this only gives us the great circles $$(e^{i\alpha},0),\quad (0,e^{i\beta})$$ respectively. 

Consider the tangent vectors $$\frac{\partial}{\partial\alpha}=\cos\theta\cdot (ie^{i\alpha},0)\quad \ \textrm{and} \ \quad \frac{\partial}{\partial\beta}=\sin\theta\cdot (0,ie^{i\beta})$$ at the point $p=(\cos\theta\cdot e^{i\alpha},\sin\theta\cdot e^{i\beta}),$ obtained through differentiation. \smallskip

We now compute the first fundamental form with respect to the Riemannian metric $g.$
\begin{eqnarray*}
E_g&=&g_p(\frac{\partial}{\partial\alpha},\frac{\partial}{\partial\alpha})\\
&=&\cos^2\theta(\lambda^2\cos^2\theta+\sin^2\theta(\mu^2\sin^2(\alpha+\beta)+\nu^2\cos^2(\alpha+\beta))),
\end{eqnarray*}
\begin{eqnarray*}
F_g&=&g_p(\frac{\partial}{\partial\alpha},\frac{\partial}{\partial\beta})\\
&=&\cos^2\theta\sin^2\theta(
\lambda^2-(\mu^2\sin^2(\alpha+\beta)+\nu^2\cos(\alpha+\beta)))\\
&=&-E_g+\lambda^2\cos^2\theta\\
&=&-G_g+\lambda^2\sin^2\theta,
\end{eqnarray*}
\begin{eqnarray*}
G_g&=&g_p(\frac{\partial}{\partial\beta},\frac{\partial}{\partial\beta})\\
&=&\sin^2\theta(\lambda^2\sin^2\theta+\cos^2\theta(\mu^2\sin^2(\alpha+\beta)+\nu^2\cos(\alpha+\beta)))\\
&=&E_g+\lambda^2(\sin^2\theta-\cos^2\theta).
\end{eqnarray*}

We now employ the Gram-Schmidt process to find an orthonormal basis $V_1,V_2$ for the tangent space $T_pT_\theta^2$ of the torus $T_\theta^2$ at $p$: 
$$
V_1=\frac{1}{\sqrt{E_g}}\cdot \frac{\partial}{\partial\alpha}.$$

\begin{eqnarray*}
V_2^{'}&=&E_g\cdot \frac{\partial}{\partial\beta}-F_g\cdot \frac{\partial}{\partial\alpha}.
\end{eqnarray*}

\begin{eqnarray*}
V_2&=&\frac{V_2^{'}}{\sqrt{g_p(V_2^{'},V_2^{'})}}\\
&=&\frac{-F_g}{\sqrt{E_g(E_g G_g-F_g^2)}}\cdot\frac{\partial}{\partial\alpha}+\frac{\sqrt{E_g}}{\sqrt{E_g G_g-F_g^2}}\cdot \frac{\partial}{\partial\beta}.
\end{eqnarray*}

We have now obtained an orthonormal basis $$V_1=f_1\cdot\frac{\partial}{\partial\alpha},\quad V_2=f_2\cdot\frac{\partial}{\partial\alpha}+f_3\cdot\frac{\partial}{\partial\beta},$$
where we define
\begin{eqnarray*}
    f_1&=&\frac{1}{\sqrt{E_g}},\\
    f_2&=&\frac{-F_g}{\sqrt{E_g(E_gG_g-F_g^2)}},\\
    f_3&=&\frac{\sqrt{E_g}}{\sqrt{E_gG_g-F_g^2}}
\end{eqnarray*}
as functions of $\alpha$ and $\beta.$
The basis $V_1,V_2$ can now be used in determining the trace of the second fundamental form $B$. Note that by Proposition 6.22 in \cite{Gud-Rie}, $B$ is tensorial in both arguments. (In the following, we index $B$ with $g$ to be able to later compare it to the Lorentzian case.)
\begin{eqnarray*}
\trace B_g&=&B_g(V_1,V_1)+B_g(V_2,V_2)\\
&=&B_g(f_1\cdot \frac{\partial}{\partial\alpha},f_1\cdot \frac{\partial}{\partial\alpha})\\&&\quad+B_g(f_2\cdot\frac{\partial}{\partial\alpha}+f_3\cdot\frac{\partial}{\partial\beta},f_2\cdot\frac{\partial}{\partial\alpha}+f_3\cdot\frac{\partial}{\partial\beta})\\
&=&(f_1^2+f_2^2)\cdot B_g(\frac{\partial}{\partial\alpha},\frac{\partial}{\partial\alpha})+2f_2f_3\cdot B_g(\frac{\partial}{\partial\alpha},\frac{\partial}{\partial\beta})\\
&&\quad +f_3^2\cdot B_g(\frac{\partial}{\partial\beta},\frac{\partial}{\partial\beta})\\
&=&(f_1^2+f_2^2)\cdot B_g(\frac{\partial}{\partial\alpha},\frac{\partial}{\partial\alpha})+f_3^2\cdot B_g(\frac{\partial}{\partial\beta},\frac{\partial}{\partial\beta}).\\
\end{eqnarray*}
In the last step we used that the mixed derivatives
$$\frac{\partial^2}{\partial\alpha\partial\beta},\frac{\partial^2}{\partial\beta\partial\alpha}$$ vanish, which implies that
$$B_g(\frac{\partial}{\partial\alpha},\frac{\partial}{\partial\beta})$$ does as well. 
We now simplify the separate terms of the sum. Instead of projecting $$B_g(\frac{\partial}{\partial\alpha},\frac{\partial}{\partial\alpha})$$ onto a unit normal vector $V_3$ of the torus, which we could find through another step of the Gram-Schmidt process, we use an orthogonal decomposition. This is justified since the normal vector $V_3$ of the torus, the unit vectors $V_1,V_2$ spanning the tangent plane $T_pT_\theta^2$ and $N_p$ form an orthonormal basis of $T_p\cn^2,$ i.e.
$$T_p\cn^2=T_pS^3\oplus N_pS^3=T_pT_\theta ^2\oplus N_pT_\theta ^2\oplus N_pS^3.$$
This orthogonal decomposition holds with respect to the Berger metric $g,$ for any choice of $\lambda,\mu,\nu\in\rn^+,$ since $V_1,V_2,V_3$ depend on those parameters.
\begin{eqnarray*}
B_g(\frac{\partial}{\partial\alpha},\frac{\partial}{\partial\alpha})&=&(\nab{\frac{\partial}{\partial\alpha}}{\frac{\partial}{\partial\alpha}})^\perp\\
&=&(\nab{\frac{\partial}{\partial\alpha}}{\frac{\partial}{\partial\alpha}})-g_p(\nab{\frac{\partial}{\partial\alpha}}{\frac{\partial}{\partial\alpha}},V_1)V_1-\\
&&\quad g_p(\nab{\frac{\partial}{\partial\alpha}}{\frac{\partial}{\partial\alpha}},V_2)V_2-g_p(\nab{\frac{\partial}{\partial\alpha}}{\frac{\partial}{\partial\alpha}},N_p)\cdot N_p.
\end{eqnarray*}
We now simplify the separate terms of $B_g(\frac{\partial}{\partial\alpha},\frac{\partial}{\partial\alpha}).$
\begin{eqnarray*}
\nab{\frac{\partial}{\partial\alpha}}{\frac{\partial}{\partial\alpha}}&=&(-\cos\theta\cdot e^{i\alpha},0).
\end{eqnarray*}
\begin{eqnarray*}
g_p(\nab{\frac{\partial}{\partial\alpha}}{\frac{\partial}{\partial\alpha}},V_1)
&=&f_1\cdot g_p(\nab{\frac{\partial}{\partial\alpha}}{\frac{\partial}{\partial\alpha}},\frac{\partial}{\partial\alpha})\\
&=&f_1\cdot(\mu^2-\nu^2)\sin^2\theta\cos^2\theta\sin(\alpha+\beta)\cos(\alpha+\beta)\\
&=&f_1\cdot \xi,
\end{eqnarray*}
where for brevity we set $$\xi=(\mu^2-\nu^2)\sin^2\theta\cos^2\theta\sin(\alpha+\beta)\cos(\alpha+\beta).$$ 
\begin{eqnarray*}
g_p(\nab{\frac{\partial}{\partial\alpha}}{\frac{\partial}{\partial\alpha}},V_2)&=&f_2\cdot g_p(\nab{\frac{\partial}{\partial\alpha}}{\frac{\partial}{\partial\alpha}},\frac{\partial}{\partial\alpha})+f_3\cdot g_p(\nab{\frac{\partial}{\partial\alpha}}{\frac{\partial}{\partial\alpha}},\frac{\partial}{\partial\beta})\\
&=&(f_2-f_3)(\mu^2-\nu^2)\sin^2\theta\cos^2\theta\sin(\alpha+\beta)\cos(\alpha+\beta)\\
&=&(f_2-f_3)\cdot \xi.
\end{eqnarray*}

$$
g_p(\nab{\frac{\partial}{\partial\alpha}}{\frac{\partial}{\partial\alpha}},N_p)=-\cos^2\theta.
$$
We now plug in the above results. We denote $$\frac{\partial}{\partial\theta}=(-\sin\theta\cdot e^{i\alpha}, \cos\theta\cdot e^{i\beta}),$$ which is the vector belonging to $T_pS^3$ obtained by differentiating for $\theta.$
\begin{eqnarray*}
&&B_g(\frac{\partial}{\partial\alpha},\frac{\partial}{\partial\alpha})\\
&=&\nab{\frac{\partial}{\partial\alpha}}{\frac{\partial}{\partial\alpha}}-f_1\cdot \xi\cdot V_1-(f_2-f_3)\cdot \xi\cdot V_2+\cos^2\theta\cdot N_p\\
&=&\nab{\frac{\partial}{\partial\alpha}}{\frac{\partial}{\partial\alpha}}-\xi(f_1^2+f_2^2-f_2f_3)\cdot \frac{\partial}{\partial\alpha}-\xi(f_2f_3-f_3^2)\cdot \frac{\partial}{\partial\beta}+\cos^2\theta\cdot N_p\\
&=&\sin\theta\cos\theta\cdot\frac{\partial}{\partial\theta}-\xi\cdot\frac{G_g+F_g}{E_g G_g-F_g^2}\cdot \frac{\partial}{\partial\alpha}+\xi\cdot\frac{E_g+F_g}{E_g G_g-F_g^2}\cdot \frac{\partial}{\partial\beta}\\
&=&\sin\theta\cos\theta\cdot\frac{\partial}{\partial\theta}+\frac{\xi}{E_g G_g-F_g ^2}\cdot(\lambda^2\sin^2\theta\cdot(-i\cos\theta\cdot e^{i\alpha},0)+\lambda^2\cos^2\theta\cdot(0,i\sin\theta\cdot e^{i\beta}))\\
&=&\sin\theta\cos\theta\cdot\frac{\partial}{\partial\theta}+\frac{\xi}{E_g G_g-F_g ^2}\cdot i\sin\theta\cos\theta\cdot\frac{\partial}{\partial\theta}\\
&=&\sin\theta\cos\theta\cdot\frac{\partial}{\partial\theta}\cdot\frac{\mu^2\sin^2(\alpha+\beta)+\nu^2\cos^2(\alpha+\beta)+i(\mu^2-\nu^2)\sin(\alpha+\beta)\cos(\alpha+\beta)}{\mu^2\sin^2(\alpha+\beta)+\nu^2\cos^2(\alpha+\beta)}\\
&=&\sin\theta\cos\theta\cdot\frac{\partial}{\partial\theta}\cdot\frac{(i\mu^2\sin(\alpha+\beta)+\nu^2\cos(\alpha+\beta))\cdot e^{-i(\alpha+\beta)}}{\mu^2\sin^2(\alpha+\beta)+\nu^2\cos^2(\alpha+\beta)}\\
&=&\sin\theta\cos\theta\cdot\frac{i\mu^2\sin(\alpha+\beta)+\nu^2\cos(\alpha+\beta)}{\mu^2\sin^2(\alpha+\beta)+\nu^2\cos^2(\alpha+\beta)}\cdot(-\sin\theta\cdot e^{-i\beta},\cos\theta\cdot e^{-i\alpha})\\
&=&\sin\theta\cos\theta\cdot\frac{\mu^2\sin(\alpha+\beta)}{\mu^2\sin^2(\alpha+\beta)+\nu^2\cos^2(\alpha+\beta)}\cdot p\cdot (0,i)\\
&&+\sin\theta\cos\theta\cdot\frac{\nu^2\cos(\alpha+\beta)}{\mu^2\sin^2(\alpha+\beta)+\nu^2\cos^2(\alpha+\beta)}\cdot p\cdot (0,1)\\
&=&\sin\theta\cos\theta\cdot\frac{\mu^2\sin(\alpha+\beta)}{\mu^2\sin^2(\alpha+\beta)+\nu^2\cos^2(\alpha+\beta)}\cdot Z_p\\
&&-\sin\theta\cos\theta\cdot\frac{\nu^2\cos(\alpha+\beta)}{\mu^2\sin^2(\alpha+\beta)+\nu^2\cos^2(\alpha+\beta)}\cdot Y_p.
\end{eqnarray*}

We now proceed with $B_g(\frac{\partial}{\partial\beta},\frac{\partial}{\partial\beta})$ the same way we did with $B_g(\frac{\partial}{\partial\alpha},\frac{\partial}{\partial\alpha}).$
\begin{eqnarray*}
B_g(\frac{\partial}{\partial\beta},\frac{\partial}{\partial\beta})&=&(\nab{\frac{\partial}{\partial\beta}}{\frac{\partial}{\partial\beta}})^\perp\\
&=&(\nab{\frac{\partial}{\partial\beta}}{\frac{\partial}{\partial\beta}})-g_p(\nab{\frac{\partial}{\partial\beta}}{\frac{\partial}{\partial\beta}},V_1)V_1-\\
&&g_p(\nab{\frac{\partial}{\partial\beta}}{\frac{\partial}{\partial\beta}},V_2)V_2-g_p(\nab{\frac{\partial}{\partial\beta}}{\frac{\partial}{\partial\beta}},N_p)\cdot N_p.
\end{eqnarray*}

\begin{eqnarray*}
(\nab{\frac{\partial}{\partial\beta}}{\frac{\partial}{\partial\beta}})&=&(0,-\sin\theta\cdot e^{i\beta}).
\end{eqnarray*}

\begin{eqnarray*}
g_p(\nab{\frac{\partial}{\partial\beta}}{\frac{\partial}{\partial\beta}},V_1)
&=&f_1\cdot g_p(\nab{\frac{\partial}{\partial\beta}}{\frac{\partial}{\partial\beta}},\frac{\partial}{\partial\alpha})\\
&=&f_1\cdot(-\mu^2+\nu^2)\sin^2\theta\cos^2\theta\sin(\alpha+\beta)\cos(\alpha+\beta)\\
&=&-f_1\cdot \xi.
\end{eqnarray*}
\begin{eqnarray*}
g_p(\nab{\frac{\partial}{\partial\beta}}{\frac{\partial}{\partial\beta}},V_2)&=&f_2\cdot g_p(\nab{\frac{\partial}{\partial\beta}}{\frac{\partial}{\partial\beta}},\frac{\partial}{\partial\alpha})+f_3\cdot g_p(\nab{\frac{\partial}{\partial\beta}}{\frac{\partial}{\partial\beta}},\frac{\partial}{\partial\beta})\\
&=&-(f_2-f_3)(\mu^2-\nu^2)\sin^2\theta\cos^2\theta\sin(\alpha+\beta)\cos(\alpha+\beta)\\
&=&-(f_2-f_3)\cdot \xi.
\end{eqnarray*}
\begin{eqnarray*}
g_p(\nab{\frac{\partial}{\partial\beta}}{\frac{\partial}{\partial\beta}},N_p)&=&-\sin^2\theta.
\end{eqnarray*}
Finally we obtain
\begin{eqnarray*}
B_g(\frac{\partial}{\partial\beta},\frac{\partial}{\partial\beta})&=&(\nab{\frac{\partial}{\partial\beta}}{\frac{\partial}{\partial\beta}})+f_1\cdot\xi\cdot V_1+(f_2-f_3)\cdot \xi\cdot V_2+\sin^2\theta\cdot N_p\\
&=&-B_g(\frac{\partial}{\partial\alpha},\frac{\partial}{\partial\alpha}).
\end{eqnarray*}
We now compile this information, to get the following:
\begin{eqnarray*}
&&\trace B_g\\
&=&(f_1^2+f_2^2)\cdot B_g(\frac{\partial}{\partial\alpha},\frac{\partial}{\partial\alpha})+f_3^2\cdot B_g(\frac{\partial}{\partial\beta},\frac{\partial}{\partial\beta})\\
&=&(f_1^2+f_2^2-f_3^2)\cdot B_g(\frac{\partial}{\partial\alpha},\frac{\partial}{\partial\alpha})\\
&=&\frac{G_g-E_g}{E_g G_g-F_g^2}\cdot B_g(\frac{\partial}{\partial\alpha},\frac{\partial}{\partial\alpha})\\
&=&\frac{(\sin^2\theta-\cos^2\theta)}{\cos\theta\sin\theta(\mu^2\sin^2(\alpha+\beta)+\nu^2\cos^2(\alpha+\beta))^2}\\
&&\quad\cdot(-\nu^2\cos(\alpha+\beta)\cdot Y_p+\mu^2\sin(\alpha+\beta)\cdot Z_p)\\
&=&\frac{-2}{\tan(2\theta)\cdot(\mu^2\sin^2(\alpha+\beta)+\nu^2\cos^2(\alpha+\beta))^2}\cdot(-\nu^2\cos(\alpha+\beta)\cdot Y_p+\mu^2\sin(\alpha+\beta)\cdot Z_p).\\
\end{eqnarray*}
The mean curvature vector $H_g$ is now given by \begin{eqnarray*}
H_g&=&\tfrac{1}{2}\cdot \trace B_g\\
&=&\frac{(\nu^2\cos(\alpha+\beta)\cdot Y_p-\mu^2\sin(\alpha+\beta)\cdot Z_p)}{\tan(2\theta)\cdot(\mu^2\sin^2(\alpha+\beta)+\nu^2\cos^2(\alpha+\beta))^2}.
\end{eqnarray*}
In particular, \begin{eqnarray*} \|H_g\|
&=&\sqrt{g_p(\tfrac{1}{2}\cdot\trace B_g, \tfrac{1}{2}\cdot\trace B_g)}\\
&=&\frac{\mu\nu}{|\tan(2\theta)|\cdot(\mu^2\sin^2(\alpha+\beta)+\nu^2\cos^2(\alpha+\beta))^{\tfrac{3}{2}}}.
\end{eqnarray*}
It is evident that $\|H_g\|$ does not depend on $\lambda,$ since $H_g$ is a linear combination of $Y_p$ and $Z_p.$ Note that for all $\mu,\nu,$ $\theta=\frac{\pi}{4}$ is the only choice within $(0,\frac{\pi}{2})$ for which 
$\cos(2\theta)=0$ and thus $\frac{1}{\tan(2\theta)}=0.$
In this case, we have
$$H_g\equiv0.$$
Thus the mean curvature vanishes identically if and only if $\theta=\frac{\pi}{4}.$ This renders the torus $T_{\pi/4}^2$ a minimal submanifold of $S^3$.
\begin{theorem}
Assume that for the metric $g$ on $S^3,$ as given in (\ref{metricg}), it holds that $\mu=\nu.$ Let $\mu$ be given. Then for every non-negative real number $C$ there exists a torus $T_\theta^2, $ belonging to the family of tori described above, such that its mean curvature is constant and satisfies $$\|H_g\|\equiv C.$$
For $C=0,$ this is unique, for $C>0,$ there exist two such tori.
\end{theorem}
\begin{proof}
We use the computations from above. 
In the case that $\mu=\nu,$ the above identity yields
\begin{eqnarray*}
\|H_g\|&=&
\frac{1}{|\tan(2\theta)|\cdot\mu}.
\end{eqnarray*}
Clearly, $\|H_g\|$ does not depend on $\alpha$ and $\beta,$ thus the mean curvature is constant along the torus. 
We then note that for $\theta$ within our interval of choice, $$\frac{1}{\tan(2\theta)}=\frac{\cos(2\theta)}{\sin(2\theta)}$$ assumes all values on the real line. There is only one solution for $\frac{1}{\tan(2\theta)}=0,$ which is $\theta=\frac{\pi}{4},$ as discussed above.\smallskip

For $C>0,$ we just have to solve $${|\tan(2\theta)|=\frac{1}{C\mu}},$$ which has two solutions:
one solution where $$\tan(2\theta)<0,\quad \textrm{i.e.}\quad \theta\in(\frac{\pi}{4},\frac{\pi}{2}),$$ and another solution where $$\tan(2\theta)>0,\quad \textrm{i.e.}\quad \theta\in(0,\frac{\pi}{4}).$$
\end{proof}

\begin{remark}
If we plug in $\theta=\frac{\pi}{4},$ we obtain the Clifford torus $$T_{{\pi/4}}=(\tfrac{1}{\sqrt2}\cdot e^{i\alpha},\tfrac{1}{\sqrt2}\cdot e^{i\beta}).$$
Thus, our result does not contradict the Lawson conjecture \cite{Lawson}. We however only investigated a particular family of tori, whereas the counterexamples to the Lawson conjecture given by Torralbo are undoloid-type surfaces. For a more detailed classification, we refer the reader to Theorem 1 and Remark 3 in \cite{TorralboRot}.
\end{remark}
\begin{remark}
Our result confirms the findings of \cite{deLima}. For this, set $\mu=\nu=1$ and $\lambda=\tau.$ While we chose a parametrisation of the torus depending on the $\theta\in(0,\frac{\pi}{2}),$ de Lima et al. define the torus depending on a parameter $r,$ namely for $z,w\in S^3$ $$T_r=(r\cdot z, \sqrt{1-r^2}\cdot w), \ \ r\in(0,1).$$ 
By setting $r=\sin\theta$ it can easily be   seen that this parametrisation is equivalent to ours.
\end{remark}

\chapter{The Generalised Lorentzian Berger Spheres $(S^3,h)$}
We now compute the Levi-Civita connection and the sectional curvatures of $S^3$, but with respect to the Lorentzian metric $h.$ \smallskip
\begin{definition}\label{metrich}
We equip $S^3$ with a family of left-invariant Lorentzian metrics $$\{h:\smo (TS^3)\otimes
\smo (TS^3)\mapsto \smo(S^3)\,|\,\lambda \ ,\mu \ ,\nu\in\rn^+\},$$ such that the restriction of $h$ to a point $p$ is given by
\begin{eqnarray*}
&&h_p(A,B)\\
&=&-\lambda^2\cdot\ip{p^{-1}A}{(i,0)}\ip{p^{-1}B}{(i,0)}+\mu^2\cdot\ip{p^{-1} A}{(0,-1)}\ip{p^{-1}B}{(0,-1)}\\ &&\quad+\nu^2\cdot\ip{p^{-1} A}{(0,i)}\ip{p^{-1} B}{(0,i)}+\ip{p^{-1} A}{(1,0)}\ip{p^{-1} B}{(1,0)},\end{eqnarray*}
where $A,B\in T_pS^3.$
An orthonormal frame for the tangent bundle $TS^3$ with respect to $h$ is given by $$
\{X_p=\lambda^{-1}\cdot p\cdot(i,0), \ Y_p=\mu^{-1}\cdot p\cdot(0,-1), \ Z_p=\nu^{-1}\cdot p\cdot(0,i)\}.
$$
whereas the normal bundle $NS^3$ is spanned by $N_p=p\cdot (1,0)$. Note that this coincides with the frames for $TS^3$ and $ NS^3$ of $(S^3,g)$ as given in Definition \ref{basis}.
\end{definition}

\begin{lemma}\label{LC5}
The Levi-Civita connection $\nabla$ of the Lorentzian Lie group $(S^3,h)$ satisfies 
$$\nab XX=0,\quad \nab XY=\frac{\lambda^2+\mu^2+\nu^2}{\lambda\mu\nu}\cdot Z,\quad \nab XZ=-\frac{\lambda^2+\mu^2+\nu^2}{\lambda\mu\nu}\cdot Y,$$
$$\nab YX=\frac{\lambda^2+\mu^2-\nu^2}{\lambda\mu\nu}\cdot Z,\quad \nab YY=0,\quad \nab YZ=\frac{\lambda^2+\mu^2-\nu^2}{\lambda\mu\nu}\cdot X,$$
$$\nab ZX=\frac{-\lambda^2+\mu^2-\nu^2}{\lambda\mu\nu}\cdot Y,\quad \nab ZY=\frac{-\lambda^2+\mu^2-\nu^2}{\lambda\mu\nu}\cdot X,\quad \nab ZZ=0.$$
\end{lemma}

\begin{proof}
The arguments are the same as in the Riemannian case, except taking the following Lorentzian formula into account
$$\nab AB=-h(\nab AB,X)X+h(\nab AB,Y)Y+h(\nab AB,Z)Z.$$
We see that compared to the Riemannian case, the signs of every $\lambda^2$ and every $X$ differ. 
\end{proof}

\begin{proposition}
The sectional curvatures of the Lorentzian Lie group $(S^3,h)$ satisfy 
\begin{eqnarray*}
K(X,Y)&=&h(R(X,Y)Y,X)=\frac{(\lambda^2+\mu^2-\nu^2)^2+4\nu^2(\mu^2-\nu^2)}{(\lambda\mu\nu)^2}\\
K(X,Z)&=&h(R(X,Z)Z,X)=\frac{(\lambda^2-\mu^2+\nu^2)^2-4\mu^2(\mu^2-\nu^2)}{(\lambda\mu\nu)^2}\\ 
K(Y,Z)&=&h(R(Y,Z)Z,Y)=\frac{(\lambda^2+\mu^2+\nu^2)^2+2(\lambda^4-\mu^4-\nu^4)}{(\lambda\mu\nu)^2}.
\end{eqnarray*}
\end{proposition}
\begin{proof}
The result is obtained using the formula for the Riemann curvature operator as given in (\ref{RCO}) and the Levi-Civita connection we computed in Lemma \ref{LC5}.
\end{proof}
\begin{remark}
In particular, 
\begin{equation}\label{KXYh}
K(X,Y)\leq0 \end{equation} if $\nu>\mu$ and \begin{equation}\label{KXYh2}
0<\lambda^2\leq2\nu\sqrt{\nu^2-\mu^2}-\mu^2+\nu^2,\end{equation}
with equality in \ref{KXYh} only if we have equality in \ref{KXYh2}.\smallskip
\begin{equation}\label{KXZh}
K(X,Z)\leq0\end{equation} if $\nu<\mu$ and \begin{equation}\label{KXZh2}
0<\lambda^2\leq2\mu\sqrt{\mu^2-\nu^2}-\nu^2+\mu^2,\end{equation} with equality in \ref{KXZh} only if we have equality in \ref{KXZh2}.
\begin{equation}\label{KYZh}
K(Y,Z)\leq0\end{equation} if \begin{equation}\label{KYZh2}
0<\lambda^2\leq\frac{1}{3}\cdot(2\sqrt{\mu^4-\mu^2\nu^2+\nu^4}-\mu^2-\nu^2),\end{equation}
with equality in \ref{KYZh} only if we have equality in \ref{KYZh2}.
\end{remark}
\chapter{CMC Tori in the Lorentzian $(S^3,h)$}
In this chapter we determine for which parameters we obtain CMC or minimal tori in $(S^3,h).$ We proceed similarly as in the Riemannian case. As it turns out, we obtain the same result.\smallskip

Consider again $$\frac{\partial}{\partial\alpha}, \ \frac{\partial}{\partial\beta}\in T_pS^3$$ as given in Chapter 4. We notice that in the first fundamental form  only the sign of $\lambda^2$ changes.
\begin{eqnarray*}
E_h&=&h_p(\frac{\partial}{\partial\alpha},\frac{\partial}{\partial\alpha})\\
&=&\cos^2\theta(-\lambda^2\cos^2\theta+\sin^2\theta(\mu^2\sin^2(\alpha+\beta)+\nu^2\cos^2(\alpha+\beta))),\\
F_h&=&h_p(\frac{\partial}{\partial\alpha},\frac{\partial}{\partial\beta})\\
&=&-\cos^2\theta\sin^2\theta(
\lambda^2+\mu^2\sin^2(\alpha+\beta)+\nu^2\cos(\alpha+\beta))\\
&=&-E_h-\lambda^2\cos^2\theta\\
&=&-G_h-\lambda^2\sin^2\theta,\\
G_h&=&h_p(\frac{\partial}{\partial\beta},\frac{\partial}{\partial\beta})\\
&=&\sin^2\theta(-\lambda^2\sin^2\theta+\cos^2\theta(\mu^2\sin^2(\alpha+\beta)+\nu^2\cos(\alpha+\beta)))\\
&=&E_h-\lambda^2(\sin^2\theta-\cos^2\theta).
\end{eqnarray*}

Through the Gram-Schmidt process we obtain an orthonormal basis for $T_pT_\theta^2:$
\begin{eqnarray*}
V_1&=&f_1\cdot \frac{\partial}{\partial\alpha},\\
V_2&=&f_2\cdot\frac{\partial}{\partial\alpha}+f_3\cdot\frac{\partial}{\partial\beta},\\
\end{eqnarray*}
where 
\begin{eqnarray*}
f_1&=&\frac{1}{\sqrt{E_h}},\\
f_2&=&\frac{-F_h}{\sqrt{E_h(E_h G_h-F_h^2)}},\\
f_3&=&\frac{\sqrt{E_h}}{\sqrt{E_h G_h-F_h^2}}.
\end{eqnarray*}
Note that the functions $f_1,f_2,f_3$ differ from the functions of the same name in Chapter 4 (when we were talking about the Riemannian case $(S^3,g)$), because they depend on the first fundamental form.
\begin{eqnarray*}
\trace B_h&=&B_h(V_1,V_1)+B_h(V_2,V_2)\\
&=&(f_1^2+f_2^2)\cdot B_h(\frac{\partial}{\partial\alpha},\frac{\partial}{\partial\alpha})+f_3^2\cdot B_h(\frac{\partial}{\partial\beta},\frac{\partial}{\partial\beta}).
\end{eqnarray*}

We now evaluate the operator $B.$ Naturally, the switch to the Lorentzian metric does not impact differentiation. 
\begin{eqnarray*}
B_h(\frac{\partial}{\partial\alpha},\frac{\partial}{\partial\alpha})&=&(\nab{\frac{\partial}{\partial\alpha}}{\frac{\partial}{\partial\alpha}})^\perp\\
&=&(\nab{\frac{\partial}{\partial\alpha}}{\frac{\partial}{\partial\alpha}})-h_p(\nab{\frac{\partial}{\partial\alpha}}{\frac{\partial}{\partial\alpha}},V_1)V_1-\\
&&\quad h_p(\nab{\frac{\partial}{\partial\alpha}}{\frac{\partial}{\partial\alpha}},V_2)V_2-h_p(\nab{\frac{\partial}{\partial\alpha}}{\frac{\partial}{\partial\alpha}},N_p)\cdot N_p\\
&=&\nab{\frac{\partial}{\partial\alpha}}{\frac{\partial}{\partial\alpha}}-\xi(f_1^2+f_2^2-f_2f_3)\cdot \frac{\partial}{\partial\alpha}-\xi(f_2f_3-f_3^2)\cdot \frac{\partial}{\partial\beta}+\cos^2\theta\cdot N_p\\
&=&\cos\theta\sin\theta\cdot\frac{\partial}{\partial\theta}-\xi\cdot\frac{G_h+F_h}{E_h G_h-F_h^2}\cdot \frac{\partial}{\partial\alpha}+\xi\cdot\frac{E_h+F_h}{E_h G_h-F_h^2}\cdot \frac{\partial}{\partial\beta},
\end{eqnarray*}
where, as in Chapter 4, $$\xi=(\mu^2-\nu^2)\sin^2\theta\cos^2\theta\sin(\alpha+\beta)\cos(\alpha+\beta)$$ and $$\frac{\partial}{\partial\theta}=(-\sin\theta\cdot e^{i\alpha}, \cos\theta\cdot e^{i\beta}).$$
Similarly, we find that 
\begin{eqnarray*}
B_h(\frac{\partial}{\partial\beta},\frac{\partial}{\partial\beta})&=&(\nab{\frac{\partial}{\partial\beta}}{\frac{\partial}{\partial\beta}})^\perp\\
&=&(\nab{\frac{\partial}{\partial\beta}}{\frac{\partial}{\partial\beta}})+\xi(f_1^2+f_2^2-f_2f_3)\cdot \frac{\partial}{\partial\alpha}+\xi(f_2f_3-f_3^2)\cdot \frac{\partial}{\partial\beta}+\sin^2\theta\cdot N_p\\
&=&-\cos\theta\sin\theta\cdot\frac{\partial}{\partial\theta}+\xi\cdot\frac{G_h+F_h}{E_h G_h-F_h^2}\cdot \frac{\partial}{\partial\alpha}-\xi\cdot\frac{E_h+F_h}{E_h G_h-F_h^2}\cdot \frac{\partial}{\partial\beta}\\
&=&-B_h(\frac{\partial}{\partial\alpha},\frac{\partial}{\partial\alpha}).
\end{eqnarray*}
Thus, $$\trace B_h=(f_1^2+f_2^2-f_3^2)\cdot B_h(\frac{\partial}{\partial\alpha},\frac{\partial}{\partial\alpha}).$$

We now look at the coefficients of $\frac{\partial}{\partial\alpha}$ and $\frac{\partial}{\partial\beta}$ in $B_h(\frac{\partial}{\partial\alpha},\frac{\partial}{\partial\alpha})$ and note that the sign changes in the first fundamental form (compared to $E_g, F_g,G_g$ in Chapter 4) cancel each other out. Indeed,
\begin{eqnarray*}
\frac{G_h+F_h}{E_h G_h-F_h^2}&=&\frac{-\lambda^2\sin^2\theta}{-\lambda^2\cos^2\theta\sin^2\theta(\mu^2\sin^2(\alpha+\beta)+\nu^2\cos^2(\alpha+\beta))}\\
&=&\frac{G_g+F_g}{E_g G_g-F_g^2},\\
\frac{E_h+F_h}{E_h G_h-F_h^2}&=&\frac{-\lambda^2\cos^2\theta}{-\lambda^2\cos^2\theta\sin^2\theta(\mu^2\sin^2(\alpha+\beta)+\nu^2\cos^2(\alpha+\beta))}\\
&=&\frac{E_g+F_g}{E_g G_g-F_g^2}.\\
\end{eqnarray*}
Consequently, $$B_h(\frac{\partial}{\partial\alpha},\frac{\partial}{\partial\alpha})=B_g(\frac{\partial}{\partial\alpha},\frac{\partial}{\partial\alpha}),$$ 
$$B_h(\frac{\partial}{\partial\beta},\frac{\partial}{\partial\beta})=B_g(\frac{\partial}{\partial\beta},\frac{\partial}{\partial\beta}).$$
Naturally, in the Lorentzian case we also have that
$$
B_h(\frac{\partial}{\partial\alpha},\frac{\partial}{\partial\beta})=B_g(\frac{\partial}{\partial\alpha},\frac{\partial}{\partial\beta})=0.
$$
Further,
\begin{eqnarray*}
    f_1^2+f_2^2-f_3^2&=&\frac{G_h-E_h}{E_h G_h-F_h^2}\\
    &=&\frac{-\lambda^2(\sin^2\theta-\cos^2\theta)}{-\lambda^2\cos^2\theta\sin^2\theta(\mu^2\sin^2(\alpha+\beta)+\nu^2\cos^2(\alpha+\beta))}\\
&=&\frac{G_g-E_g}{E_g G_g-F_g^2}.
\end{eqnarray*}
This means that 
\begin{eqnarray*}
\trace B_h
&=&\frac{G_h-E_h}{E_h G_h-F_h^2}\cdot B_h(\frac{\partial}{\partial\alpha},\frac{\partial}{\partial\alpha})\\
&=&\frac{G_g-E_g}{E_g G_g-F_g^2}\cdot B_g(\frac{\partial}{\partial\alpha},\frac{\partial}{\partial\alpha})\\
&=&\trace B_g.
\end{eqnarray*}
Since $\trace B_g$ is orthogonal to $X_p,$ $$h_p(\trace B_g,\trace B_g)=g_p(\trace B_g,\trace B_g).$$
This implies that  $\|H_h\|$ in the Lorentzian case will coincide with its Riemannian counterpart $\|H_g\|.$ We obtain the same result: namely that the torus $T_\theta^2$ is minimal if and only if $\theta=\frac{\pi}{4},$ and that the mean curvature is constant if $\mu=\nu$ or $\theta=\frac{\pi}{4}.$
\begin{remark}
To compute the mean curvature in the Maple program that can be found in the Appendix, only set $\epsilon_1=-1.$ Other semi-Riemannian cases can easily be checked by changing the other $\epsilon s.$
\end{remark}

\chapter{The Generalised Riemannian Dual Space $(\Sigma^3,g)$}
We first consider a subset $\Sigma^3$ of $\cn^2$ and the diffeomorphism $\Psi:\Sigma^3\mapsto G\subset\SLC 2.$ Further, we show that $G$ is diffeomorphic to $\SLR 2.$ We then compute the Levi-Civita connection and sectional curvatures. \smallskip

Consider the subset $\Sigma^3=\{(z,w)\in\cn^2\,|\,|z|^2-|w|^2=1\}$ of $\cn^2$. Let $\Psi:\Sigma^3\to\GLC 2$ be the diffeomorphism $$\Psi:(z,w)\mapsto\begin{pmatrix}
    z&\bar w\\
    w&\bar z
\end{pmatrix}.$$
For any $(z,w)\in\Sigma^3$, $$\det(\Psi(z,w))=|z|^2-|w|^2=1.$$
This shows that the image $G$ of $\Sigma^3$ under $\Psi$ is contained in the complex special linear group $$\mathbf{SL}_2(\cn)=\{X\in\cn^{2\times2}\,|\,\det X=1\}.$$\smallskip

The standard matrix multiplication $*$ on $\mathbf{SL}_2(\cn)$ induces a natural group structure $\cdot$ on $\Sigma^3.$ 
The multiplication $\cdot$ on $\Sigma^3$ is thus defined by $$(z_1,w_1)\cdot(z_2,w_2)=(z_1z_2+\bar w_1 w_2,w_1z_2+\bar z_1 w_2).$$
The multiplicative inverse of an element $p=(z,w)\in\Sigma^3$ is given by $$p^{-1}=(\bar z,-w).$$
Assume now that $$p=\begin{pmatrix}
    z_1&\bar w_1\\
    w_1&\bar z_1
\end{pmatrix}, \ \  q=\begin{pmatrix}
    z_2&\bar w_2\\
    w_2&\bar z_2,
\end{pmatrix} \in G=\Psi(\Sigma^3).$$ Then $$p*q^{-1}=\begin{pmatrix}
    z_1\bar z_2-\bar w_1w_2&-z_1\bar w_2+\bar w_1z_2\\
    w_1\bar z_2-\bar z_1w_2&-w_1\bar w_2+\bar z_1z_2
\end{pmatrix}.$$
Thus $p*q^{-1}\in G,$ and consequently $G$ is a subgroup of $\mathbf{SL}_2(\cn).$ By proposition 2.35 in \cite{Gud-Rie}, $(G,*)$ is a Lie group.\\

Further we have that $(G,*)$ is isomorphic to $(\mathbf{SL}_2(\rn),*).$ To show this, we write $z=x+iy, \  w=a+ib,$ and let $\hat\Psi:G\to\mathbf{SL}_2(\rn)$ be the diffeomorphism given by$$\hat\Psi:
\begin{pmatrix}
x+iy & a-ib\\
a+ib & x-iy
\end{pmatrix}
\mapsto
\begin{pmatrix}
x-a & b-y\\
b+y & x+a
\end{pmatrix}.$$
Indeed, the determinant of the image of $\hat\Psi$ is 1: \begin{eqnarray*}\det \begin{pmatrix}
x-a & b-y\\
b+y & x+a
\end{pmatrix}&=&(x+a)(x-a)-(b+y)(b-y)\\
&=&x^2+y^2-a^2-b^2\\
&=&|z|^2-|w|^2\\
&=&1.\end{eqnarray*}\\
The standard scalar product of $\SLR 2$ is given by $$\ip{A}{B}=\tfrac{1}{2}\cdot\trace{A^t\cdot B}.$$
The Lie algebra $\slr 2$ is generated by
$$\{\begin{pmatrix} 0&-1\\1&0\end{pmatrix}, \begin{pmatrix}
1&0\\0&-1
\end{pmatrix}, \begin{pmatrix}
0&1\\1&0
\end{pmatrix}\}.$$
From this it can be seen that an orthonormal basis for the tangent space of $\Sigma^3$ at the neutral element $e$ with respect to the scalar product on $\rn^4,$ $$\ip{(z_1,w_1)}{(z_2,w_2)}=\Re(\bar z_1 z_2+\bar w_1 w_2).$$ is given by $$\{(i,0),\ (0,-1),\ (0,i)\}.$$
We now equip $\Sigma^3$ with a family of left invariant Riemannian Berger metrics 
$$\{g:\smo (T\Sigma^3)\otimes 
\smo (T\Sigma^3)\to \smo(\Sigma^3)\,|\,\lambda \ ,\mu \ ,\nu\in\rn^+\},$$ 
such that for each $p\in \Sigma^3$ the restriction of $g$ to $p$ is a real scalar product
$$g_p:T_p \Sigma^3\otimes 
T_p \Sigma^3\to\rn,$$ given by 
\begin{eqnarray*}
g_p(A,B)&=&\lambda^2\cdot\ip{p^{-1} A}{(i,0)}\ip{p^{-1} B}{(i,0)}+\mu^2\cdot\ip{p^{-1} A}{(0,-1)}\ip{p^{-1}B}{(0,-1)}\\
& &\quad +\,\nu^2\cdot\ip{p^{-1} A}{(0,i)}\ip{p^{-1} B}{(0,i)}+\ip{p^{-1} A}{(1,0)}\ip{p^{-1}B}{(1,0)},
\end{eqnarray*}
for $A,B\in T_p\Sigma^3.$
With respect to the Berger metric $g$ above, an orthonormal frame for the tangent bundle $T\Sigma^3$ of $\Sigma^3$ is now given by
$$
\{X_p=\lambda^{-1}\cdot p\cdot(i,0), \quad Y_p=\mu^{-1}\cdot p\cdot(0,-1), \quad Z_p=\nu^{-1}\cdot p\cdot(0,i)\},
$$
the normal bundle $N\Sigma^3$ is spanned by $N_p=p\cdot(1,0).$ Matrix multiplication on $\slr 2$ yields the bracket relations
$$\lb XY=2\, \lambda^{-1}\mu^{-1}\nu\, Z,\quad \lb ZX=2\, \lambda^{-1}\mu\nu^{-1}\, Y,\quad \lb YZ=-2\, \lambda\mu^{-1}\nu^{-1}\, X.$$

\begin{lemma}
The Levi-Civita connection $\nabla$ of the Riemannian Lie group $(\Sigma^3,g)$ satisfies 
$$\nab XX=0,\quad \nab XY=\frac{\lambda^2+\mu^2+\nu^2}{\lambda\mu\nu}\cdot Z,\quad \nab XZ=-\frac{\lambda^2+\mu^2+\nu^2}{\lambda\mu\nu}\cdot Y,$$
$$\nab YX=\frac{\lambda^2+\mu^2-\nu^2}{\lambda\mu\nu}\cdot Z,\quad \nab YY=0,\quad \nab YZ=-\frac{\lambda^2+\mu^2-\nu^2}{\lambda\mu\nu}\cdot X,$$
$$\nab ZX=\frac{-\lambda^2+\mu^2-\nu^2}{\lambda\mu\nu}\cdot Y,\quad \nab ZY=-\frac{-\lambda^2+\mu^2-\nu^2}{\lambda\mu\nu}\cdot X,\quad \nab ZZ=0.$$
\end{lemma}

\begin{proof}
This is a standard calculation employing the Koszul formula (\ref{Koszul})
for left-invariant vector fields and the Lie bracket relations.
\end{proof}

\begin{proposition}
The sectional curvatures of the Riemannian Lie group $(\Sigma^3,g)$ satisfy 
\begin{eqnarray*}
K(X,Y)&=&g(R(X,Y)Y,X)=\frac{(\lambda^2+\mu^2-\nu^2)^2+4\nu^2(\mu^2-\nu^2)}{(\lambda\mu\nu)^2},\\
K(X,Z)&=&g(R(X,Z)Z,X)=\frac{(\lambda^2-\mu^2+\nu^2)^2-4\mu^2(\mu^2-\nu^2)}{(\lambda\mu\nu)^2},\\ 
K(Y,Z)&=&g(R(Y,Z)Z,Y)=-\frac{(\lambda^2+\mu^2+\nu^2)^2+2(\lambda^4-\mu^2-\nu^2)}{(\lambda\mu\nu)^2}.
\end{eqnarray*}
\end{proposition}
\begin{proof}
We computed the Riemann curvature operator as given in (\ref{RCO}). It simplifies the computation immensely if we compare the Levi-Civita connection on $(\Sigma^3,g)$ and $(S^3,h)$ and note that the signs in front of every $X$ differ (see Lemma \ref{LC5}). There is another sign change for the Lie bracket $\lb YZ$. Finally, for the Riemann curvature operator on $(\Sigma^3,g)$
we yield -1 times of what we obtained on $(S^3,h).$
\end{proof}
\begin{remark}
\begin{equation} \label{eq:KXY7}
K(X,Y)\leq0\end{equation} if $\nu>\mu$ and \begin{equation}\label{eq:KXY72}
0<\lambda^2\leq2\nu\sqrt{\nu^2-\mu^2}-\mu^2+\nu^2,\end{equation}
with equality in (\ref{eq:KXY7}) only if we have equality in (\ref{eq:KXY72}).
\begin{equation}\label{eq:KXZ7}
K(X,Z)\leq0\end{equation} if $\nu<\mu$ and \begin{equation}\label{eq:KXZ72}
0<\lambda^2\leq2\mu\sqrt{\mu^2-\nu^2}-\nu^2+\mu^2,\end{equation}
with equality in (\ref{eq:KXZ7}) only if we have equality in (\ref{eq:KXZ72})
\begin{equation}\label{eq:KYZ7}
K(Y,Z)\leq0\end{equation} if \begin{equation}\label{eq:KYZ72}
\lambda^2\geq\frac{1}{3}(2\sqrt{\mu^4-\mu^2\nu^2+\nu^4}-\mu^2-\nu^2),\end{equation}
with equality in (\ref{eq:KYZ7}) only if we have equality in (\ref{eq:KYZ72}).
\end{remark}



\chapter{CMC surfaces in the Riemannian $(\Sigma^3,g)$}
In this chapter we consider a family of tori in $\Sigma^3$ and determine for which parameters they have CMC or are minimal. Using a similar strategy as in Chapter 4, we then compute the mean curvature.  \smallskip

For a fixed $\theta>0,$ we parametrise a two dimensional submanifold with the map $\mathcal{F}_\theta:\rn^2\to\Sigma^3$ defined by $$\mathcal{F}_\theta: (\alpha,\beta)\to (\cosh\theta\cdot e^{i\alpha},\sinh\theta\cdot e^{i\beta}).$$
This yields a family of tori $U_\theta^2,$ which we equip with the metric $g$ as given in the previous chapter. To justify our choice of interval for $\theta$, note that setting $\theta=0$ would again only give us the circle $(e^{i\alpha},0).$\smallskip

Consider now the tangent vectors at a point $p\in U_\theta ^2$ $$\frac{\partial}{\partial\alpha}=\cosh\theta\cdot(ie^{i\alpha},0), \quad \frac{\partial}{\partial\beta}=\sinh\theta\cdot(0,ie^{i\beta}).$$
The first fundamental form is given as follows:
\begin{eqnarray*}
E_g&=&g_p(\frac{\partial}{\partial\alpha},\frac{\partial}{\partial\alpha})\\
&=&\cosh^2\theta(\lambda^2\cosh^2\theta+\sinh^2\theta(\mu^2\sin^2(\alpha+\beta)+\nu^2\cos^2(\alpha+\beta))),
\end{eqnarray*}
\begin{eqnarray*}
F_g&=&g_p(\frac{\partial}{\partial\alpha},\frac{\partial}{\partial\beta})\\
&=&-\cosh^2\theta\sinh^2\theta(
\lambda^2+\mu^2\sin^2(\alpha+\beta)+\nu^2\cos^2(\alpha+\beta))\\
&=&-E_g+\lambda^2\cosh^2\theta\\
&=&-G_g-\lambda^2\sinh^2\theta,
\end{eqnarray*}
\begin{eqnarray*}
G_g&=&g_p(\frac{\partial}{\partial\beta},\frac{\partial}{\partial\beta})\\
&=&\sinh^2\theta(\lambda^2\sinh^2\theta+\cosh^2\theta(\mu^2\sin^2(\alpha+\beta)+\nu^2\cos(\alpha+\beta)))\\
&=&E_g-\lambda^2(\cosh^2\theta+\sinh^2\theta).
\end{eqnarray*}

We now employ the Gram-Schmidt process to find an orthonormal basis $V_1,V_2$ for the tangent space $T_p\Sigma^3$ of the torus $U_\theta^2$ at a point $p$:
$$
V_1=\frac{1}{\sqrt{E_g}}\cdot\frac{\partial}{\partial\alpha},$$
\begin{eqnarray*}
{V_2}^{'}&=&E_g\cdot \frac{\partial}{\partial\beta}-F_g\cdot \frac{\partial}{\partial\alpha},
\end{eqnarray*}
\begin{eqnarray*}
V_2&=&\frac{{V_2}^{'}}{\sqrt{g_p({V_2}^{'},{V_2}^{'})}}\\
&=&\frac{-F_g}{\sqrt{E_g (E_g G_g-F_g^2)}}\cdot\frac{\partial}{\partial\alpha}+\frac{\sqrt{E_g}}{\sqrt{E_g G_g-F_g^2}}\cdot \frac{\partial}{\partial\beta}.
\end{eqnarray*}
For simplicity, we write 
$$V_1=f_1\cdot \frac{\partial}{\partial\alpha},\quad V_2=f_2\cdot\frac{\partial}{\partial\alpha}+f_3\cdot\frac{\partial}{\partial\beta},$$
where 
\begin{eqnarray*}
f_1&=&\frac{1}{\sqrt{E_g}},\\
f_2&=&\frac{-F_g}{\sqrt{E_g (E_g G_g-F_g^2)}},\\
f_3&=&\frac{\sqrt{E_g}}{\sqrt{E_g G_g-F_g^2}}
\end{eqnarray*}
are functions of $\alpha$ and $\beta.$ The basis $V_1,V_2$ can now be used in determining the trace of $B_g$. 
\begin{eqnarray*}
\trace B_g&=&B_g(V_1,V_1)+B_g(V_2,V_2)\\
&=&B(f_1\cdot \frac{\partial}{\partial\alpha},f_1\cdot \frac{\partial}{\partial\alpha})\\&&\quad+B_g(f_2\cdot\frac{\partial}{\partial\alpha}+f_3\cdot\frac{\partial}{\partial\beta},f_2\cdot\frac{\partial}{\partial\alpha}+f_3\cdot\frac{\partial}{\partial\beta})\\
&=&(f_1^2+f_2^2)\cdot B_g(\frac{\partial}{\partial\alpha},\frac{\partial}{\partial\alpha})+2f_2f_3\cdot B_g(\frac{\partial}{\partial\alpha},\frac{\partial}{\partial\beta})\\
&&\quad +f_3^2\cdot B_g(\frac{\partial}{\partial\beta},\frac{\partial}{\partial\beta})\\
&=&(f_1^2+f_2^2)\cdot B_g(\frac{\partial}{\partial\alpha},\frac{\partial}{\partial\alpha})+f_3^2\cdot B_g(\frac{\partial}{\partial\beta},\frac{\partial}{\partial\beta}).\\
\end{eqnarray*}
Here we used that $B$ is tensorial.
Also note that since the mixed derivatives
$$\frac{\partial^2}{\partial\alpha\partial\beta}\quad \textrm{and}\quad\frac{\partial^2}{\partial\beta\partial\alpha}\quad \textrm{vanish, so does}\quad
B(\frac{\partial}{\partial\alpha},\frac{\partial}{\partial\beta}).$$\\
We now simplify the sum. As in both $S^3$ cases, we use an orthogonal decomposition instead of projecting $$B_g(\frac{\partial}{\partial\alpha},\frac{\partial}{\partial\alpha})$$ onto a unit normal vector $V_3$ of $U_\theta^2$. 
\begin{eqnarray*}
B_g(\frac{\partial}{\partial\alpha},\frac{\partial}{\partial\alpha})&=&(\nab{\frac{\partial}{\partial\alpha}}{\frac{\partial}{\partial\alpha}})^\perp\\
&=&\nab{\frac{\partial}{\partial\alpha}}{\frac{\partial}{\partial\alpha}}-g_p(\nab{\frac{\partial}{\partial\alpha}}{\frac{\partial}{\partial\alpha}},V_1)V_1\\
&&\quad -g_p(\nab{\frac{\partial}{\partial\alpha}}{\frac{\partial}{\partial\alpha}},V_2)V_2-g_p(\nab{\frac{\partial}{\partial\alpha}}{\frac{\partial}{\partial\alpha}},N_p)\cdot N_p.
\end{eqnarray*}
We proceed by simplifying the separate terms of $B_g(\frac{\partial}{\partial\alpha},\frac{\partial}{\partial\alpha}).$
\begin{eqnarray*}
\nab{\frac{\partial}{\partial\alpha}}{\frac{\partial}{\partial\alpha}}&=&(-\cosh\theta\cdot e^{i\alpha},0),
\end{eqnarray*}

\begin{eqnarray*}
g_p(\nab{\frac{\partial}{\partial\alpha}}{\frac{\partial}{\partial\alpha}},V_1)
&=&f_1\cdot g_p(\nab{\frac{\partial}{\partial\alpha}}{\frac{\partial}{\partial\alpha}},\frac{\partial}{\partial\alpha})\\
&=&f_1\cdot(\mu^2-\nu^2)\sinh^2\theta\cosh^2\theta\sin(\alpha+\beta)\cos(\alpha+\beta)\\
&=&f_1\cdot \zeta.
\end{eqnarray*}
For brevity we set $$\zeta=(\mu^2-\nu^2)\sinh^2\theta\cosh^2\theta\sin(\alpha+\beta)\cos(\alpha+\beta).$$ 
\begin{eqnarray*}
g_p(\nab{\frac{\partial}{\partial\alpha}}{\frac{\partial}{\partial\alpha}},V_2)&=&f_2\cdot g_p(\nab{\frac{\partial}{\partial\alpha}}{\frac{\partial}{\partial\alpha}},\frac{\partial}{\partial\alpha})+f_3\cdot g_p(\nab{\frac{\partial}{\partial\alpha}}{\frac{\partial}{\partial\alpha}},\frac{\partial}{\partial\beta})\\
&=&(f_2-f_3)(\mu^2-\nu^2)\sinh^2\theta\cosh^2\theta\sin(\alpha+\beta)\cos(\alpha+\beta)\\
&=&(f_2-f_3)\cdot \zeta,
\end{eqnarray*}
and finally
$$
g_p(\nab{\frac{\partial}{\partial\alpha}}{\frac{\partial}{\partial\alpha}},N_p)=-\cosh^2\theta.
$$
We now plug in the above results. In the following chain of equations we write $\frac{\partial}{\partial\theta}=(\sinh\theta\cdot e^{i\alpha},\cosh\theta\cdot e^{i\beta}).$
\begin{eqnarray*}
&&B_g(\frac{\partial}{\partial\alpha},\frac{\partial}{\partial\alpha})\\
&=&\nab{\frac{\partial}{\partial\alpha}}{\frac{\partial}{\partial\alpha}}-f_1\cdot \zeta\cdot V_1-(f_2-f_3)\cdot \zeta\cdot V_2+\cosh^2\theta\cdot N_p\\
&=&\nab{\frac{\partial}{\partial\alpha}}{\frac{\partial}{\partial\alpha}}-\zeta(f_1^2+f_2^2-f_2f_3)\cdot \frac{\partial}{\partial\alpha}-\zeta(f_2f_3-f_3^2)\cdot \frac{\partial}{\partial\beta}+\cosh^2\theta\cdot N_p\\
&=&\cosh\theta\sinh\theta\cdot\frac{\partial}{\partial\theta}-\zeta\cdot\frac{G_g+F_g}{E_g G_g-F_g^2}\cdot \frac{\partial}{\partial\alpha}+\zeta\cdot\frac{E_g+F_g}{E_g G_g-F_g^2}\cdot \frac{\partial}{\partial\beta}\\
&=&\cosh\theta\sinh\theta\cdot\frac{\partial}{\partial\theta}+\frac{(\mu^2-\nu^2)i\sin(\alpha+\beta)\cos(\alpha+\beta)\cosh\theta\sinh\theta}{\mu^2\sin^2(\alpha+\beta)+\nu^2\cos^2(\alpha+\beta)}\cdot\frac{\partial}{\partial\theta}\\
&=&\cosh\theta\sinh\theta\cdot\frac{\partial}{\partial\theta}\cdot\frac{(\mu^2i\sin(\alpha+\beta)+\nu^2\cos(\alpha+\beta))\cdot e^{-i(\alpha+\beta)}}{\mu^2\sin^2(\alpha+\beta)+\nu^2\cos^2(\alpha+\beta)}\\
&=&\cosh\theta\sinh\theta\cdot\frac{\mu^2i\sin(\alpha+\beta)+\nu^2\cos(\alpha+\beta)}{\mu^2\sin^2(\alpha+\beta)+\nu^2\cos^2(\alpha+\beta)}\cdot(\sinh\theta\cdot e^{-i\beta},\cosh\theta\cdot e^{-i\alpha})\\
&=&\cosh\theta\sinh\theta\cdot\frac{\mu^2\sin(\alpha+\beta)}{\mu^2\sin^2(\alpha+\beta)+\nu^2\cos^2(\alpha+\beta)}\cdot \nu Z_p\\
&&-\cosh\theta\sinh\theta\cdot\frac{\nu^2\cos(\alpha+\beta)}{\mu^2\sin^2(\alpha+\beta)+\nu^2\cos^2(\alpha+\beta)}\cdot\mu Y_p.\\
\end{eqnarray*}

We now proceed with $B_g(\frac{\partial}{\partial\beta},\frac{\partial}{\partial\beta})$ the same way we did with $B_g(\frac{\partial}{\partial\alpha},\frac{\partial}{\partial\alpha}).$
\begin{eqnarray*}
B_g(\frac{\partial}{\partial\beta},\frac{\partial}{\partial\beta})&=&(\nab{\frac{\partial}{\partial\beta}}{\frac{\partial}{\partial\beta}})^\perp\\
&=&(\nab{\frac{\partial}{\partial\beta}}{\frac{\partial}{\partial\beta}})-g_p(\nab{\frac{\partial}{\partial\beta}}{\frac{\partial}{\partial\beta}},V_1)V_1-\\
&&g_p(\nab{\frac{\partial}{\partial\beta}}{\frac{\partial}{\partial\beta}},V_2)V_2-g_p(\nab{\frac{\partial}{\partial\beta}}{\frac{\partial}{\partial\beta}},N_p)\cdot N_p.
\end{eqnarray*}
We now compute the different terms of $B_g(\frac{\partial}{\partial\beta},\frac{\partial}{\partial\beta}).$
\begin{eqnarray*}
\nab{\frac{\partial}{\partial\beta}}{\frac{\partial}{\partial\beta}}&=&(0,-\sinh\theta\cdot e^{i\beta}).
\end{eqnarray*}
\begin{eqnarray*}
g_p(\nab{\frac{\partial}{\partial\beta}}{\frac{\partial}{\partial\beta}},V_1)
&=&f_1\cdot g_p(\nab{\frac{\partial}{\partial\beta}}{\frac{\partial}{\partial\beta}},\frac{\partial}{\partial\alpha})\\
&=&f_1\cdot(-\mu^2+\nu^2)\sinh^2\theta\cosh^2\theta\sin(\alpha+\beta)\cos(\alpha+\beta)\\
&=&-f_1\cdot \zeta,
\end{eqnarray*}
where $\zeta$ is defined as above.
\begin{eqnarray*}
g_p(\nab{\frac{\partial}{\partial\beta}}{\frac{\partial}{\partial\beta}},V_2)
&=&f_2\cdot g_p(\nab{\frac{\partial}{\partial\beta}}{\frac{\partial}{\partial\beta}},\frac{\partial}{\partial\alpha})+f_3\cdot g_p(\nab{\frac{\partial}{\partial\beta}}{\frac{\partial}{\partial\beta}},\frac{\partial}{\partial\beta})\\
&=&(-f_2+f_3)\cdot \zeta,
\end{eqnarray*}
\begin{eqnarray*}
g_p(\nab{\frac{\partial}{\partial\beta}}{\frac{\partial}{\partial\beta}},N_p)
&=&\sinh^2\theta.
\end{eqnarray*}
From this it can be seen that
\begin{eqnarray*}
B_g(\frac{\partial}{\partial\beta},\frac{\partial}{\partial\beta})&=&\nab{\frac{\partial}{\partial\beta}}{\frac{\partial}{\partial\beta}}+\zeta(f_1^2+f_2^2-f_2f_3)\cdot \frac{\partial}{\partial\alpha}+\zeta(f_2f_3-f_3^2)\cdot \frac{\partial}{\partial\beta}\\
&&\quad -\sinh^2\theta\cdot N_p\\
&=&-\cosh\theta\sinh\theta\cdot\frac{\partial}{\partial\theta}+\zeta\cdot\frac{G_g+F_g}{E_g G_g-F_g^2}\cdot \frac{\partial}{\partial\alpha}-\zeta\cdot\frac{E_g+F_g}{E_g G_g-F_g^2}\cdot \frac{\partial}{\partial\beta}\\
&=&-B_g(\frac{\partial}{\partial\alpha},\frac{\partial}{\partial\alpha}).
\end{eqnarray*}
We can now use the above identities to determine $H_g$:
\begin{eqnarray*}
H_g&=&\frac{1}{2}\cdot\trace B_g\\
&=&\frac{1}{2}\cdot(f_1^2+f_2^2-f_3^2)\cdot B_g(\frac{\partial}{\partial\alpha},\frac{\partial}{\partial\alpha})\\
&=&\frac{1}{2}\cdot\frac{G_g-E_g}{E_g G_g-F_g^2}\cdot B_g(\frac{\partial}{\partial\alpha},\frac{\partial}{\partial\alpha})\\
&=&-\frac{\cosh^2\theta+\sinh^2\theta}{2\cosh^2\theta\sinh^2\theta(\mu^2\sin^2(\alpha+\beta)+\nu^2\cos^2(\alpha+\beta))}\cdot B_g(\frac{\partial}{\partial\alpha},\frac{\partial}{\partial\alpha})\\
&=&\frac{(\cosh^2\theta+\sinh^2\theta)\mu\nu}{2\cosh\theta\sinh\theta(\mu^2\sin^2(\alpha+\beta)+\nu^2\cos^2(\alpha+\beta))^2}\\
&&\quad\cdot(\nu\cos(\alpha+\beta)\cdot Y_p-\mu\sin(\alpha+\beta)\cdot Z_p).
\end{eqnarray*}
We finally obtain
\begin{eqnarray*}
\|H_g\|&=&\frac{(\cosh^2\theta+\sinh^2\theta)\mu\nu}{2\cosh\theta\sinh\theta(\mu^2\sin^2(\alpha+\beta)+\nu^2\cos^2(\alpha+\beta))^{\frac{3}{2}}}\\
&=&\frac{\mu\nu}{\tanh(2\theta)\cdot(\mu^2\sin^2(\alpha+\beta)+\nu^2\cos^2(\alpha+\beta))^{\frac{3}{2}}}.
\end{eqnarray*}
There exists no real $\theta$ satisfying $$\cosh^2\theta+\sinh^2\theta=0.$$
This means that none of the tori $U_\theta^2,$ as described above, is minimal. 

\begin{theorem}
Equip $\Sigma^3$ with a metric $g$ as given in the previous chapter, such that $\mu=\nu.$
Then for every real number $C>\frac{1}{\mu},$ there exists a torus $U_\theta^2$ in $(\Sigma^3,g)$ as described above, such that its mean curvature is constant and satisfies $$\|H_g\|\equiv C.$$
\end{theorem}
\begin{proof}
We use the computations from above. 
In the case that $\mu=\nu,$ the above identity yields
\begin{eqnarray*}
\|H_g\|&=&\frac{1}{|\tanh(2\theta)|\cdot\mu}.
\end{eqnarray*}
Clearly, $\|H_g\|$ does not depend on $\alpha$ and $\beta,$ thus the mean curvature is constant along all torus $U_\theta^2$ belonging to the family. Further, for all $\theta>0,$ 
$$\frac{1}{\tanh(2\theta)}>1.$$
We now only have to solve $$C=\frac{1}{\tanh(2\theta)\cdot\mu},$$
which has the unique solution
$$\theta=\frac{1}{2}\cdot\textrm{arctanh}(\frac{1}{C\mu}).$$
\end{proof}


\chapter{The Generalised Lorentzian Dual Space $(\Sigma^3,h)$} 
We equip $\Sigma^3$ with a Lorentzian metric $h$, and compute the Levi-Civita connection, as well as the sectional curvatures on $(\Sigma^3,h).$\smallskip

Let $\Sigma^3=\{(z,w)\in\cn^2\,|\,|z|^2-|w|^2=1\}$. We recall from Chapter 7 that the multiplication $\cdot$ is defined by $$(z_1,w_1)\cdot(z_2,w_2)=(z_1z_2+\bar w_1 w_2,w_1z_2+\bar z_1 w_2),$$
and $\ip{}{}$ denotes the scalar product on $\rn^4,$
 $$\ip{(z_1,w_1)}{(z_2,w_2)}=\Re(\bar z_1 z_2+\bar w_1 w_2).$$
\begin{definition}
Now equip $\Sigma^3$ with a family of left-invariant Lorentzian metrics $$\{h:\smo (T\Sigma^3)\otimes
\smo (T\Sigma^3)\to \smo(\Sigma^3)\,|\,\lambda \ ,\mu \ ,\nu\in\rn^+\},$$ such that the restriction of $h$ to a point $p\in\Sigma^3$ is given by
\begin{eqnarray*}
&&h_p(A,B)\\
&=&-\lambda^2\cdot\ip{p^{-1}A}{(i,0)}\ip{p^{-1}B}{(i,0)}+\mu^2\cdot\ip{p^{-1} A}{(0,-1)}\ip{p^{-1}B}{(0,-1)}\\ &&\quad+\nu^2\cdot\ip{p^{-1} A}{(0,i)}\ip{p^{-1} B}{(0,i)}+\ip{p^{-1} A}{(1,0)}\ip{p^{-1} B}{(1,0)},\end{eqnarray*}
where $A,B\in T_p\Sigma^3.$
An orthonormal frame of the tangent bundle $T\Sigma^3$ with respect to $h$ is given by 
$$
\{X=\lambda^{-1}\cdot p\cdot(i,0), \quad Y=\mu^{-1}\cdot p\cdot(0,-1), \quad Z=\nu^{-1}\cdot p\cdot(0,i)\}.
$$
The normal bundle $N\Sigma^3$ is spanned by $N_p=p\cdot(1,0).$
\end{definition}\smallskip

The bracket relations
$$\lb XY=2\, \lambda^{-1}\mu^{-1}\nu\, Z\quad \lb ZX=2\, \lambda^{-1}\mu\nu^{-1}\, Y,\quad \lb YZ=-2\, \lambda\mu^{-1}\nu^{-1}\, X$$ still hold as in the group $(\Sigma^3,g).$

\begin{lemma}\label{LC9}
	The Levi-Civita connection $\nabla$ of the Lorentzian Lie group $(\Sigma^3,h)$ satisfies 
	$$\nab XX=0,\quad \nab XY=\frac{-\lambda^2+\mu^2+\nu^2}{\lambda\mu\nu}\cdot Z,\quad \nab XZ=-\frac{-\lambda^2+\mu^2+\nu^2}{\lambda\mu\nu}\cdot Y,$$
$$\nab YX=\frac{-\lambda^2+\mu^2-\nu^2}{\lambda\mu\nu}\cdot Z,\quad \nab YY=0,\quad \nab YZ=-\frac{\lambda^2-\mu^2+\nu^2}{\lambda\mu\nu}\cdot X,$$
$$\nab ZX=\frac{\lambda^2+\mu^2-\nu^2}{\lambda\mu\nu}\cdot Y,\quad \nab ZY=\frac{\lambda^2+\mu^2-\nu^2}{\lambda\mu\nu}\cdot X,\quad \nab ZZ=0.$$
\end{lemma}
\begin{proof}
The arguments are the same as in the other cases. We used the bracket relations and the following Lorentzian formula 
$$\nab AB=-h(\nab AB,X)X+h(\nab AB,Y)Y+h(\nab AB,Z)Z.$$
\end{proof}

\begin{proposition}
The sectional curvature $K$ of the Lorentzian Lie group \\ $(\Sigma^3,h)$ satisfies 
\begin{eqnarray*}
K(X,Y)&=&h(R(X,Y)Y,X)=\frac{(\lambda^2-\mu^2+\nu^2)^2+4\nu^2(\mu^2-\nu^2)}{(\lambda\mu\nu)^2},\\
K(X,Z)&=&h(R(X,Z)Z,X)=\frac{(\lambda^2+\mu^2-\nu^2)^2-4\mu^2(\mu^2-\nu^2)}{(\lambda\mu\nu)^2},\\
K(Y,Z)&=&h(R(Y,Z)Z,Y)=-\frac{(\lambda^2+\mu^2+\nu^2)^2-4(\lambda^4+\mu^2\nu^2)}{(\lambda\mu\nu)^2}.\\
\end{eqnarray*}
\end{proposition}
\begin{proof}
We compute the Riemann curvature tensor as given in (\ref{RCO}) using the results of Lemma \ref{LC9}. Note that in this step we obtain -1 times our results in $(S^3,g),$ due to sign switches in the Levi-Civita connection and the Lie-bracket (compare with Proposition \ref{prop:curvatures3}).
\end{proof}
\begin{remark}
\begin{equation}\label{eq:KXY91}
    K(X,Y)\leq0
\end{equation} if $\mu<\nu$ and
\begin{equation}\label{eq:KXY92}
    0<\lambda^2\leq 2\nu\sqrt{\nu^2-\mu^2}+\mu^2-\nu^2,
\end{equation}
with equality in (\ref{eq:KXY91}) only if we have equality in  (\ref{eq:KXY92}). 
\begin{equation}\label{eq:KXZ91}
    K(X,Z)\leq0
\end{equation} if $\mu>\nu$ and
\begin{equation}\label{eq:KXZ92}
    0<\lambda^2\leq 2\mu\sqrt{\mu^2-\nu^2}+\nu^2-\mu^2,
\end{equation}
with equality in (\ref{eq:KXZ91}) only if we have equality in (\ref{eq:KXZ92}). 
\begin{equation}\label{eq:KYZ91}
    K(Y,Z)\leq0
\end{equation} if \begin{equation}\label{eq:KYZ92}
    0<\lambda^2\leq\frac{1}{3}\cdot(2\sqrt{\mu^4-\mu^2\nu^2+\nu^4}+\mu^2+\nu^2),
\end{equation}
with equality in (\ref{eq:KYZ91}) only if we have equality in (\ref{eq:KYZ92}). 
\end{remark}

\chapter{CMC surfaces in the Lorentzian $(\Sigma^3,h)$ }

In this chapter we consider a family of tori in $\Sigma^3$ and evaluate their mean curvature, this time with respect to the Lorentzian metric $h$. In particular, we determine the parameters for which they are minimal and have CMC. Due to sign cancellations, we obtain the same results for $(\Sigma^3,g)$ and $(\Sigma^3,h).$\smallskip

For a fixed $\theta>0,$
we parametrise the torus $U_\theta^2\subset\Sigma^3$ with the map $F_\theta:\rn^2\to\Sigma^3$ defined by $$F_\theta: (\alpha,\beta)\to (\cosh\theta\cdot e^{i\alpha},\sinh\theta\cdot e^{i\beta}).$$
Let $$\frac{\partial}{\partial\alpha}=\cosh\theta\cdot(ie^{i\alpha},0),\quad \frac{\partial}{\partial\beta}=\sinh\theta\cdot(0,ie^{i\beta})$$ be tangent vectors of $U_\theta^2$ at a point $p.$
The first fundamental form is given as follows:
\begin{eqnarray*}
E_h&=&h_p(\frac{\partial}{\partial\alpha},\frac{\partial}{\partial\alpha})\\
&=&\cosh^2\theta(-\lambda^2\cosh^2\theta+\sinh^2\theta(\mu^2\sin^2(\alpha+\beta)+\nu^2\cos^2(\alpha+\beta))),\end{eqnarray*}
\begin{eqnarray*}
F_h&=&h_p(\frac{\partial}{\partial\alpha},\frac{\partial}{\partial\beta})\\
&=&\cosh^2\theta\sinh^2\theta(
\lambda^2-(\mu^2\sin^2(\alpha+\beta)+\nu^2\cos(\alpha+\beta)))\\
&=&-E_h-\lambda^2\cosh^2\theta\\
&=&-G_h+\lambda^2\sinh^2\theta,
\end{eqnarray*}
\begin{eqnarray*}
G_h&=&h_p(\frac{\partial}{\partial\beta},\frac{\partial}{\partial\beta})\\
&=&\sinh^2\theta(-\lambda^2\sinh^2\theta+\cosh^2\theta(\mu^2\sin^2(\alpha+\beta)+\nu^2\cos(\alpha+\beta)))\\
&=&E_h+\lambda^2(\cosh^2\theta+\sinh^2\theta).
\end{eqnarray*}

Note that compared to the Riemannian case only the sign in front of $\lambda^2$ changes.
Through the Gram-Schmidt process we obtain an orthonormal basis $V_1,V_2$ for the tangent space $T_p U_\theta^2$: 
$$V_1=f_1\cdot \frac{\partial}{\partial\alpha},\quad V_2=f_2\cdot\frac{\partial}{\partial\alpha}+f_3\cdot\frac{\partial}{\partial\beta},$$
where \begin{eqnarray*}
f_1&=&\frac{1}{\sqrt{E_h}},\\
f_2&=&\frac{-F_h}{\sqrt{E_h (E_h G_h-F_h^2)}},\\
f_3&=&\frac{\sqrt{E_h}}{\sqrt{E_h G_h-F_h^2}}
\end{eqnarray*}
are functions of $\alpha$ and $\beta.$ Note that the functions $f_1, f_2 ,f_3$ differ from the functions of the same name in $(\Sigma^3,g),$ since they depend on the first fundamental form.
With the basis $V_1,V_2$ we can now determine $H_h$. 
\begin{eqnarray*}
\trace B_h&=&B_h(V_1,V_1)+B_h(V_2,V_2)\\
&=&(f_1^2+f_2^2)\cdot B_h(\frac{\partial}{\partial\alpha},\frac{\partial}{\partial\alpha})+f_3^2\cdot B_h(\frac{\partial}{\partial\beta},\frac{\partial}{\partial\beta}).\\
\end{eqnarray*}
Naturally, like in the Riemannian case,
$B_h(\frac{\partial}{\partial\alpha},\frac{\partial}{\partial\beta})$ vanishes. Thus we only need to simplify $B_h(\frac{\partial}{\partial\alpha},\frac{\partial}{\partial\alpha})$ and $ B_h(\frac{\partial}{\partial\beta},\frac{\partial}{\partial\beta}).$
\begin{eqnarray*}
B_h(\frac{\partial}{\partial\alpha},\frac{\partial}{\partial\alpha})&=&(\nab{\frac{\partial}{\partial\alpha}}{\frac{\partial}{\partial\alpha}})^\perp\\
&=&\nab{\frac{\partial}{\partial\alpha}}{\frac{\partial}{\partial\alpha}}-h_p(\nab{\frac{\partial}{\partial\alpha}}{\frac{\partial}{\partial\alpha}},V_1)V_1\\
&&\quad -h_p(\nab{\frac{\partial}{\partial\alpha}}{\frac{\partial}{\partial\alpha}},V_2)V_2-h_p(\nab{\frac{\partial}{\partial\alpha}}{\frac{\partial}{\partial\alpha}},N_p)\cdot N_p.
\end{eqnarray*}
We evaluate the different terms. 
\begin{eqnarray*}
h_p(\nab{\frac{\partial}{\partial\alpha}}{\frac{\partial}{\partial\alpha}},V_1)&=&
f_1\cdot h_p(\nab{\frac{\partial}{\partial\alpha}}{\frac{\partial}{\partial\alpha}},\frac{\partial}{\partial\alpha})\\
&=&f_1\cdot\zeta,\\
h_p(\nab{\frac{\partial}{\partial\alpha}}{\frac{\partial}{\partial\alpha}},V_2)&=&f_2\cdot h_p(\nab{\frac{\partial}{\partial\alpha}}{\frac{\partial}{\partial\alpha}},\frac{\partial}{\partial\alpha})+f_3\cdot h_p(\nab{\frac{\partial}{\partial\alpha}}{\frac{\partial}{\partial\alpha}},\frac{\partial}{\partial\beta})\\
&=&(f_2-f_3)\cdot\zeta,
\end{eqnarray*}
where $$\zeta=(\mu^2-\nu^2)\sinh^2\theta\cosh^2\theta\sin(\alpha+\beta)\cos(\alpha+\beta).$$
\begin{eqnarray*}
h_p(\nab{\frac{\partial}{\partial\alpha}}{\frac{\partial}{\partial\alpha}},N_p)&=&-\cosh^2\theta.
\end{eqnarray*} 
 Plugging this in, we now yield the following:
\begin{eqnarray*}
&&B_h(\frac{\partial}{\partial\alpha},\frac{\partial}{\partial\alpha})\\
&=&\nab{\frac{\partial}{\partial\alpha}}{\frac{\partial}{\partial\alpha}}-\zeta(f_1^2+f_2^2-f_2f_3)\cdot \frac{\partial}{\partial\alpha}-\zeta(f_2f_3-f_3^2)\cdot \frac{\partial}{\partial\beta}+\cosh^2\theta\cdot N_p\\
&=&\cosh\theta\sinh\theta\cdot\frac{\partial}{\partial\theta}-\zeta\cdot\frac{G_h+F_h}{E_h G_h-F_h^2}\cdot \frac{\partial}{\partial\alpha}+\zeta\cdot\frac{E_h+F_h}{E_h G_h-F_h^2}\cdot \frac{\partial}{\partial\beta},\\
\end{eqnarray*}
where $\frac{\partial}{\partial\theta}=(\sinh\theta\cdot e^{i\alpha},\cosh\theta\cdot e^{i\beta}).$\smallskip

Note now that
\begin{eqnarray*}
E_h G_h-F_h^2&=&-\lambda^2\cosh^2\theta\sinh^2\theta(\mu^2\sin^2(\alpha+\beta)+\nu^2\cos^2(\alpha+\beta))\\
&=&-(E_g G_g-F_g^2),\\
G_h+F_h&=&\lambda^2\sinh^2\theta=-(G_g+F_g),\\
E_h+F_h&=&-\lambda^2\cosh^2\theta=-(E_g+F_g).\\
\end{eqnarray*}
The sign changes (compared to the Riemannian case $(\Sigma^3,g)$) in both the numerator and denominator now cancel each other out, i.e.
\begin{eqnarray*}
B_h(\frac{\partial}{\partial\alpha},\frac{\partial}{\partial\alpha})&=&\cosh\theta\sinh\theta\cdot\frac{\partial}{\partial\theta}-\zeta\cdot\frac{G_g+F_g}{E_g G_g-F_g^2}\cdot \frac{\partial}{\partial\alpha}+\zeta\cdot\frac{E_g+F_g}{E_g G_g-F_g^2}\cdot \frac{\partial}{\partial\beta}\\
&=&B_g(\frac{\partial}{\partial\alpha},\frac{\partial}{\partial\alpha}).
\end{eqnarray*} 
Through a similar procedure, we can confirm that
\begin{eqnarray*}
B_h(\frac{\partial}{\partial\beta},\frac{\partial}{\partial\beta})&=&B_g(\frac{\partial}{\partial\beta},\frac{\partial}{\partial\beta})\\
&=&-B_h(\frac{\partial}{\partial\alpha},\frac{\partial}{\partial\alpha}).\\
\end{eqnarray*}

Finally, we can compile this information:
\begin{eqnarray*}
\trace B_h&=&(f_1^2+f_2^2-f_3^2)\cdot B_h(\frac{\partial}{\partial\alpha},\frac{\partial}{\partial\alpha})\\
&=&\frac{G_h-E_h}{E_h G_h-F_h^2}\cdot B_h(\frac{\partial}{\partial\alpha},\frac{\partial}{\partial\alpha})\\
&=&\frac{G_g-E_g}{E_g G_g-F_g^2}\cdot B_g(\frac{\partial}{\partial\alpha},\frac{\partial}{\partial\alpha})\\
&=&\trace B_g.
\end{eqnarray*}
It follows that $H_h=H_g.$
We find out that changing from the Riemannian to the Lorentzian metric did not impact the mean curvature vector, nor its norm, since $B_h(\frac{\partial}{\partial\alpha},\frac{\partial}{\partial\alpha})$ is a linear combination of the vectors $Y_p$ and $Z_p,$ which are orthogonal to $X_p.$
\begin{eqnarray*}
\|H_h\|&=&\frac{(\cosh^2\theta+\sinh^2\theta)\mu\nu}{2\cosh\theta\sinh\theta(\mu^2\sin^2(\alpha+\beta)+\nu^2\cos^2(\alpha+\beta))^{\frac{3}{2}}}\\
&=&\frac{\mu\nu}{\tanh(2\theta)\cdot(\mu^2\sin^2(\alpha+\beta)+\nu^2\cos^2(\alpha+\beta))^{\frac{3}{2}}}.
\end{eqnarray*}
Thus, we conclude that there exist no minimal tori $U_\theta^2$ in $(\Sigma^3,h)$, and that all tori $U_\theta^2$ have constant mean curvature if and only if $\mu=\nu.$

\appendix
\chapter{Computations in Maple}
\newpage
\pagenumbering{gobble}
\begin{figure}[htbp]
\hspace*{-2cm}   
    \adjustimage{scale=1,left}{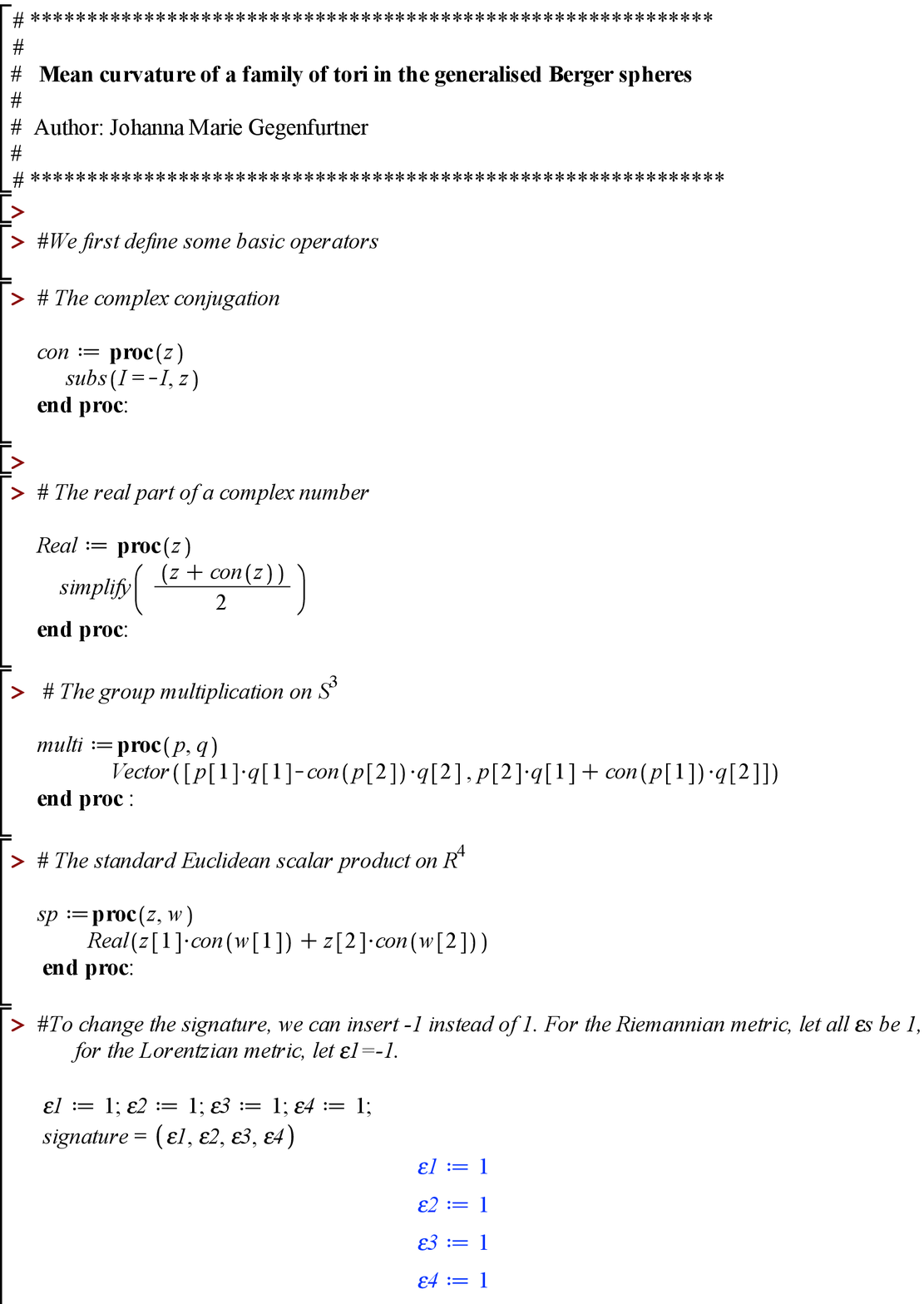}
\end{figure}
\begin{figure}[htbp]
\hspace*{-2cm}   
    \adjustimage{scale=1,left}{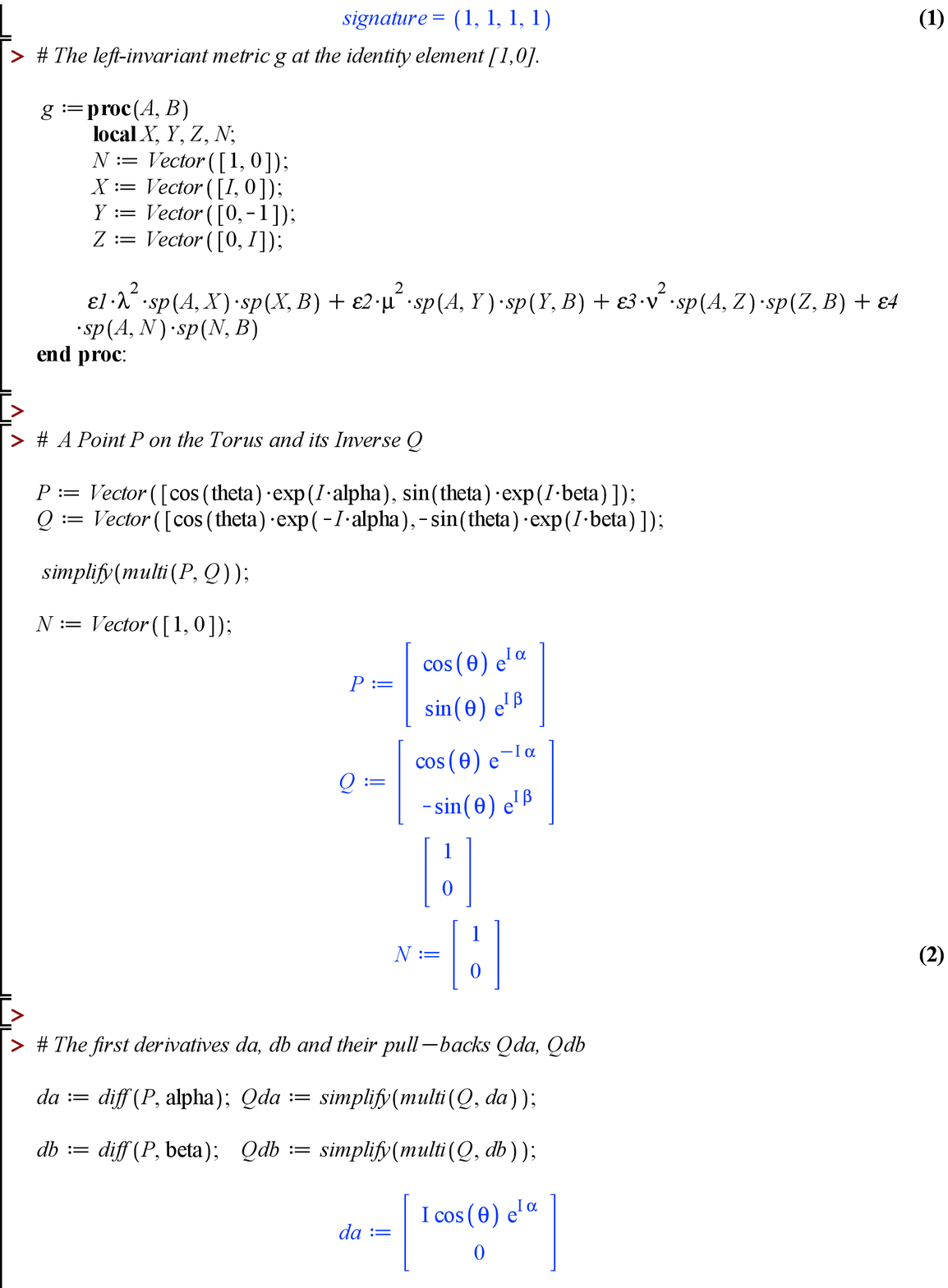}
\end{figure}
\begin{figure}[htbp]
\hspace*{-2cm}   
    \adjustimage{scale=1,left}{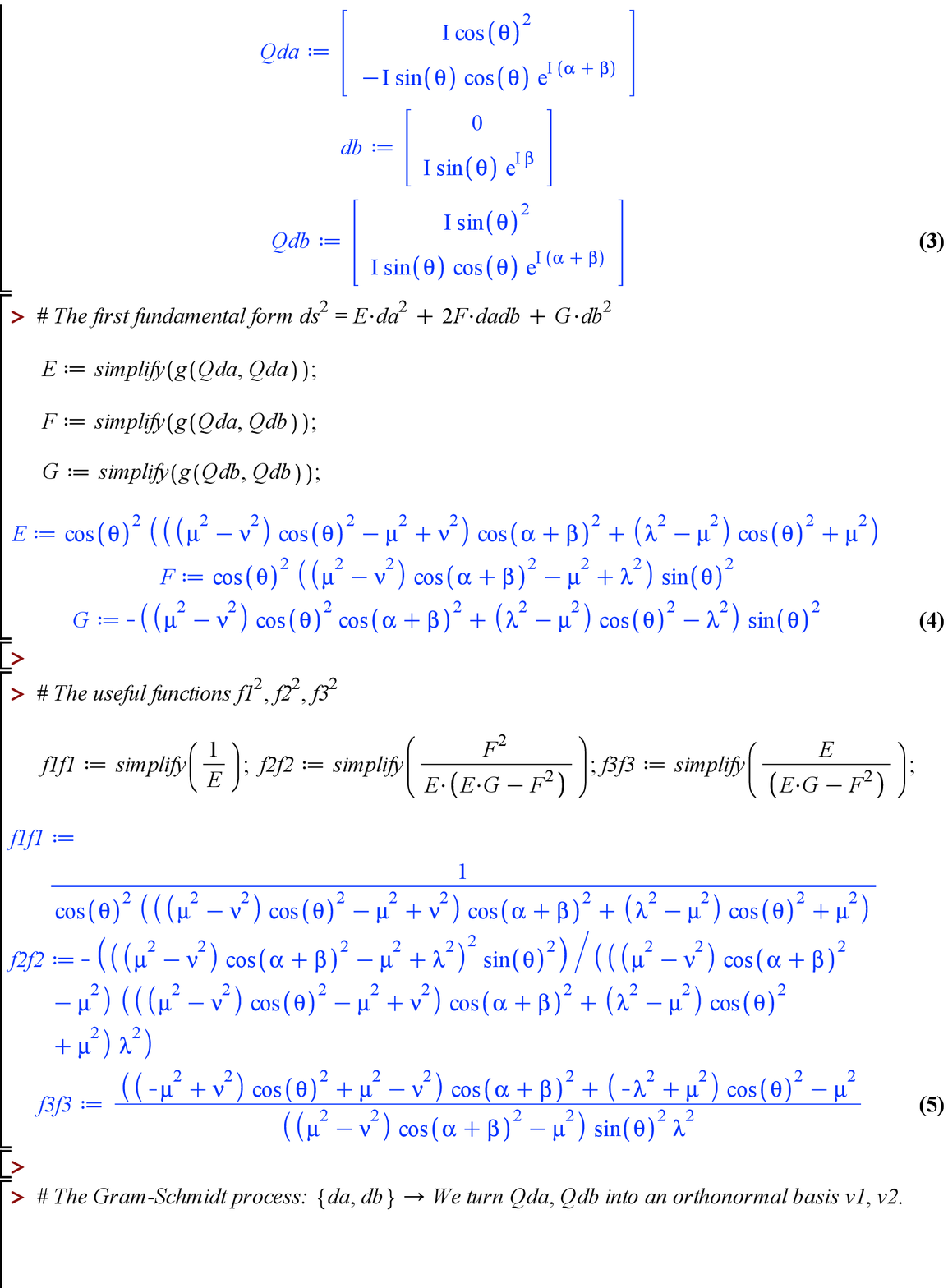}
\end{figure}
\begin{figure}[htbp]
\hspace*{-2cm}   
    \adjustimage{scale=1,left}{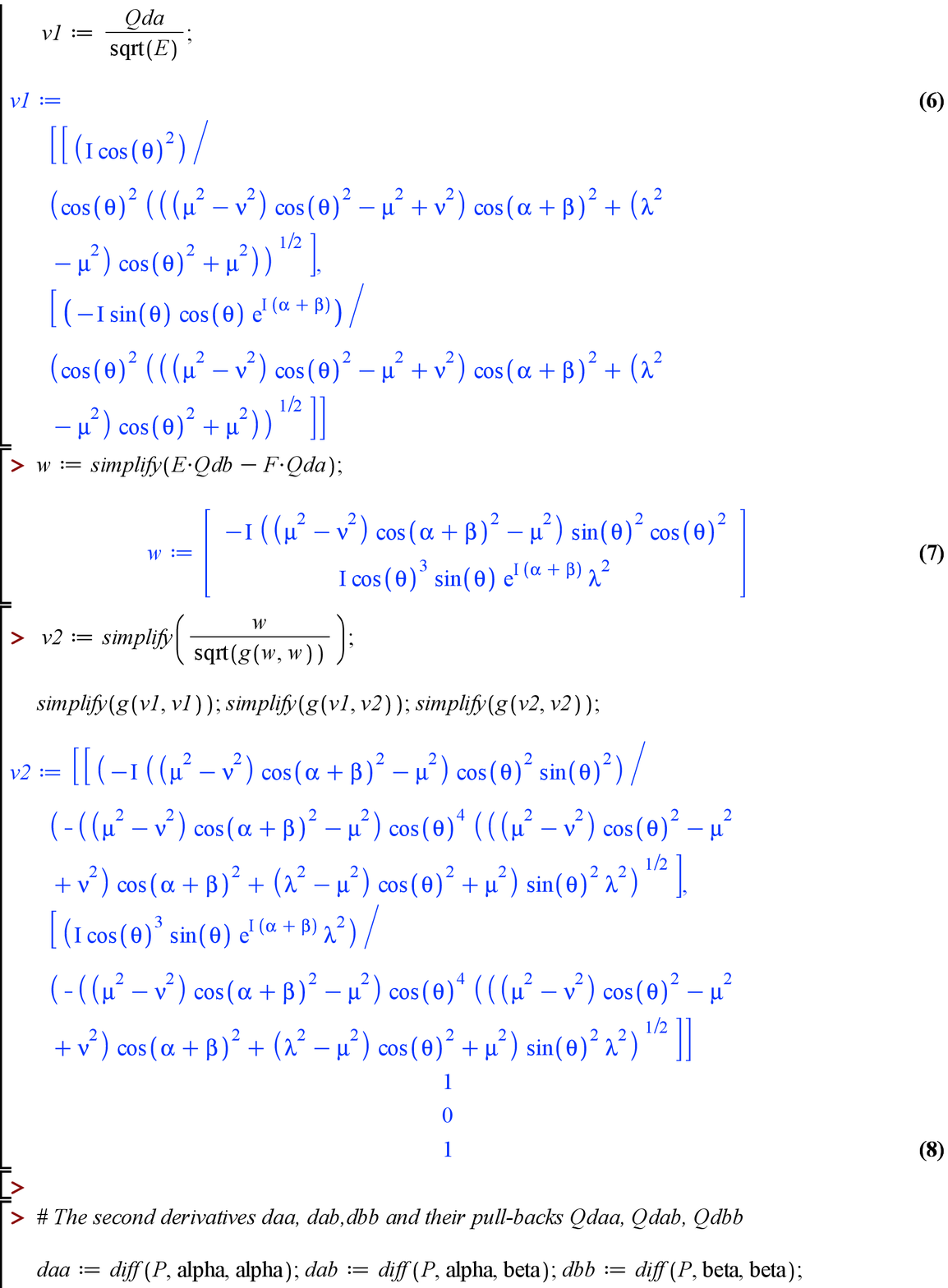}
\end{figure}
\begin{figure}[htbp]
\hspace*{-2cm}   
    \adjustimage{scale=1,left}{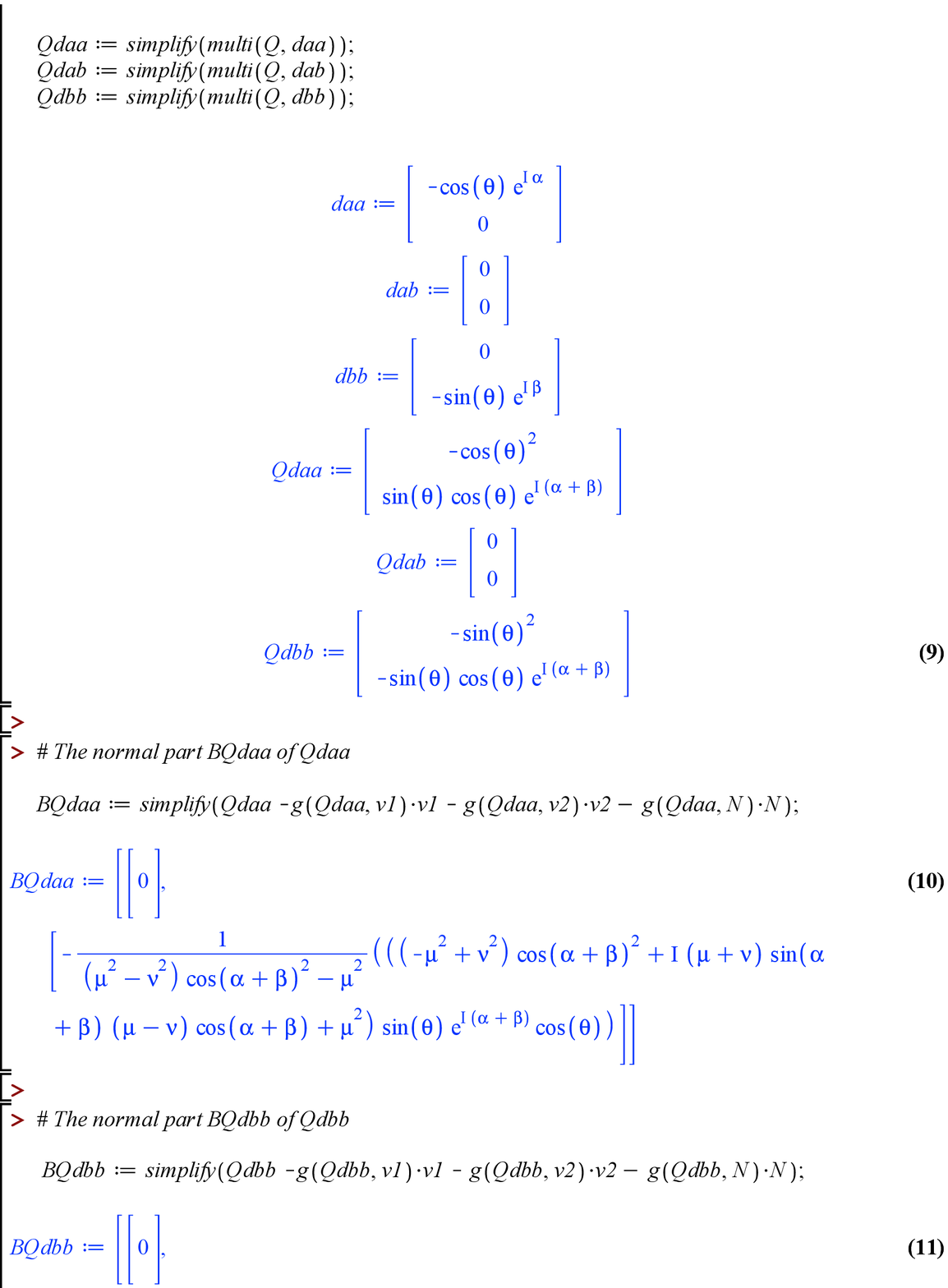}
\end{figure}
\begin{figure}[htbp]
\hspace*{-2cm}   
    \adjustimage{scale=1,left}{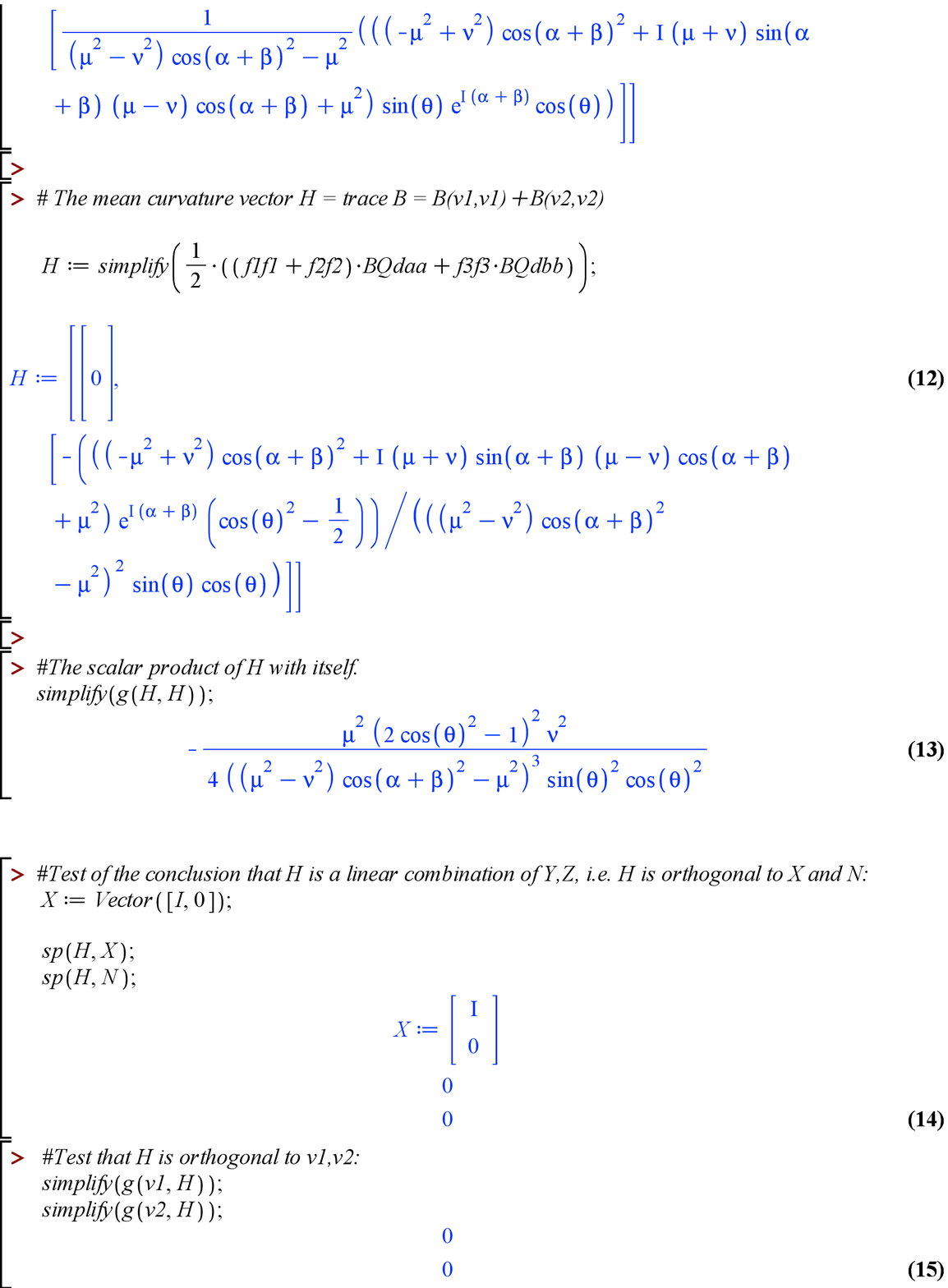}
\end{figure}
\begin{figure}[htbp]
\vspace*{-8cm}
\hspace*{-2cm}   
    \adjustimage{scale=1,left}{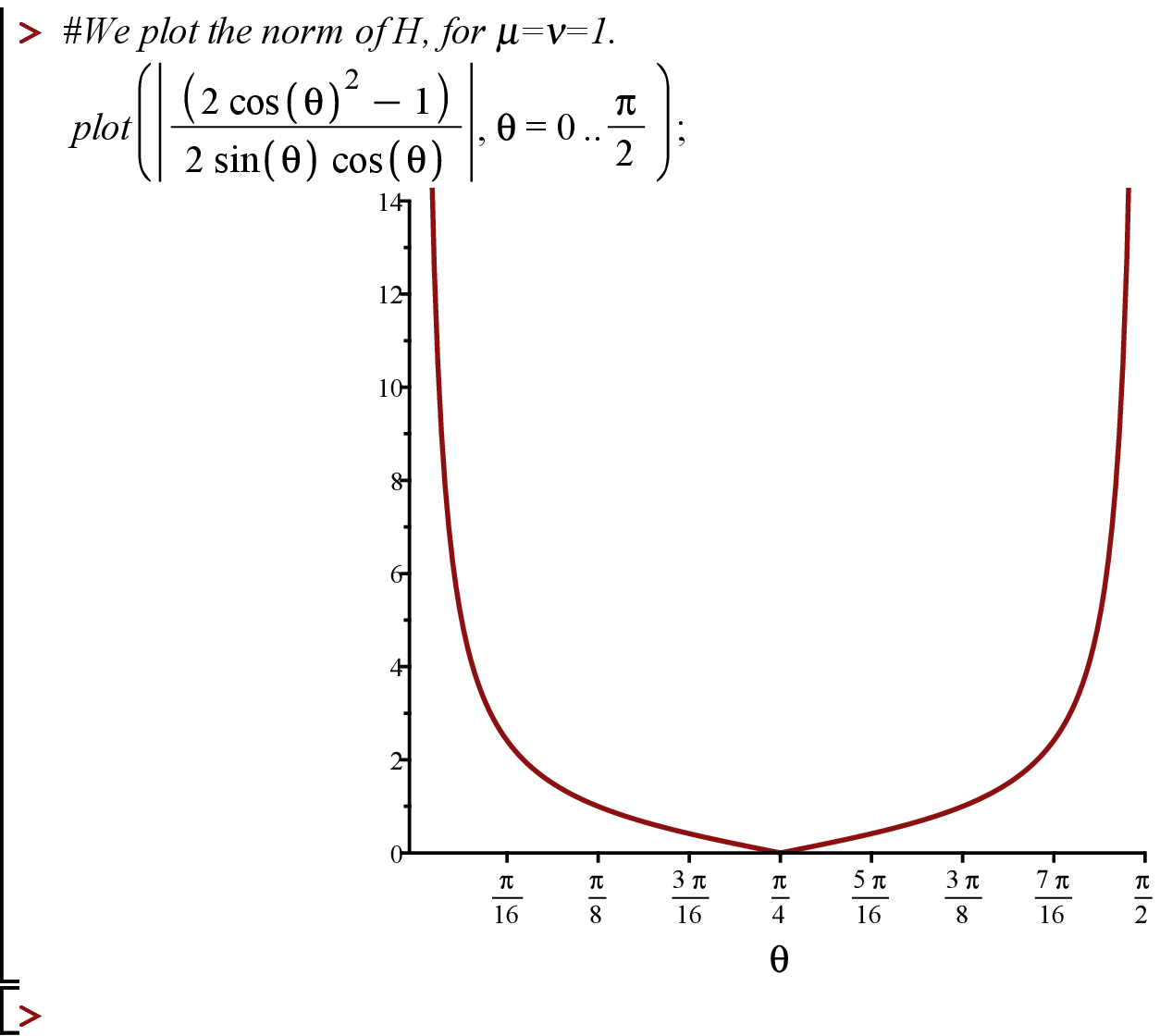}
\end{figure}

\backcover

\end{document}